\documentclass[11pt]{article}

\usepackage[breaklinks]{hyperref}
\usepackage[caption=false,font=scriptsize]{subfig}
\usepackage{amsmath, bbm}
\usepackage{amssymb}
\usepackage{graphicx}
\usepackage{setspace}
\usepackage{makecell}
\usepackage{epsfig}
\usepackage{amsfonts}
\usepackage[square,sort,comma,numbers]{natbib}
\usepackage{color}
\usepackage{mathrsfs}
\usepackage{epstopdf}
\usepackage{xcolor}
\usepackage{colortbl}
\usepackage{multirow}
\usepackage{booktabs}
\usepackage{lscape}
\usepackage{caption}
\usepackage{algorithm}
\usepackage[noend]{algpseudocode}
\usepackage{amsthm} 
\usepackage{theoremref}
\usepackage{enumerate}
\usepackage[inline]{enumitem}

\usepackage{csvsimple}
\usepackage[top=1in, bottom=1in, left=1in, right=1in]{geometry}
\usepackage{diagbox}
\usepackage{csvsimple,booktabs,array,siunitx}

\usepackage{breakurl}
\usepackage{mathtools}
\usepackage{authblk}
\usepackage{scalefnt}
\usepackage{tabstackengine}
\usepackage{float}

\usepackage{hyperref}
\usepackage{chngcntr}
\usepackage{apptools}
\AtAppendix{\counterwithin{lemma}{section}}
\AtAppendix{\counterwithin{proposition}{section}}
\AtAppendix{\counterwithin{remark}{section}}
\AtAppendix{\counterwithin{example}{section}}
\AtAppendix{\counterwithin{table}{section}}

\allowdisplaybreaks

\doublespace

\newtheorem{definition}{Definition}
\newtheorem{theorem}{Theorem}
\newtheorem{lemma}{Lemma}

\newtheorem{proposition}{Proposition}
\newtheorem{corollary}{Corollary}
\newtheorem{example}{Example}
\newtheorem{remark}{Remark}

\newcommand{\cA}{\mathcal{A}}

\newcommand{\cB}{\mathcal{B}}
\newcommand{\cC}{\mathcal{C}}

\newcommand{\cL}{\mathcal{L}}
\newcommand{\cS}{\mathcal{S}}
\newcommand{\cT}{\mathcal{T}}

\newcommand{\cW}{\mathcal{W}}
\newcommand{\cK}{\mathcal{K}}

\newcommand{\cZ}{\mathcal{Z}}
\newcommand{\forAll}{\forall\text{ }}
\newcommand{\setbar}{\text{ }|\text{ }}
\newcommand{\myProof}{\noindent\textbf{Proof: }}

\newcommand{\myf}{f}

\newcommand{\myh}{h}

\newcommand{\metric}{\zeta}
\newcommand{\glob}{\mathcal{E}}
\newcommand{\conv}{\mathrm{conv}}
\newcommand{\dknap}{D}

\providecommand{\dotminus}{\mathbin{\mathpalette\xdotminus\relax}}
\newcommand{\xdotminus}[2]{%
  \ooalign{\hidewidth$\vcenter{\hbox{$#1\dot{}$}}$\hidewidth\cr$#1-$\cr}%
}
 
\newcommand*{\myQED}{\hfill\ensuremath{\square}}

\newif\ifNoCut
\NoCutfalse
\def\cutOut#1{\ifNoCut{\color{red}#1}\fi}

\newenvironment{makeshiftResult}[1]
{
\noindent\textsc{\thnameref{#1}~\ref{#1}.}\em}
{

}

\makeatletter
\csvset{
  autobooktabularcenter/.style={
    file=#1,
    after head=\csv@pretable\begin{tabular}{*{\csv@columncount}{c}}\csv@tablehead,
    table head=\toprule\csvlinetotablerow\\\midrule,
    late after line=\\,
    table foot=\\\bottomrule,
    late after last line=\csv@tablefoot\end{tabular}\csv@posttable,
    command=\csvlinetotablerow},
}
\makeatother
\newcommand{\csvautobooktabularcenter}[2][]{\csvloop{autobooktabularcenter={#2},#1}}

\providecommand{\keywords}[1]
{
  \small	
  \textbf{\textit{Keywords---}} #1
}

\renewenvironment{abstract}
{\begin{quote}
\noindent \rule{\linewidth}{.5pt}\par{\bfseries \abstractname.}}
{\medskip\noindent \rule{\linewidth}{.5pt}
\end{quote}
}

\title{A Note on the Implications of Approximate Submodularity in Discrete Optimization}
\author[1]{Temitayo Ajayi}
\author[2]{Taewoo Lee \thanks{tlee6@uh.edu}}
\author[1]{Andrew J. Schaefer}

\affil[1]{\footnotesize Nature Source Improved Plants}
\affil[2]{\footnotesize Department of Industrial Engineering, University of Houston}
\affil[3]{\footnotesize Rice University, Department of Computational and Applied Mathematics}

\date{}

\newenvironment{repeattheorem}[1]
{
\noindent\textsc{\thnameref{#1}~\ref{#1}.}\em}
{

}

\begin{document}
\maketitle

\begin{abstract}
\small
Submodularity is a key property in discrete optimization. Submodularity has been widely used for analyzing the greedy algorithm to give performance bounds and providing insight into the construction of valid inequalities for mixed-integer programs. In recent years, researchers started to study approximate submodularity, with a primary focus on providing performance bounds for iterative approaches. In this paper, we study approximate submodularity from a different perspective in order to broaden its use cases in discrete optimization. We define metrics that quantify approximate submodularity, which we then use to derive new properties about both approximate submodularity preservation and the well-known Lov\'{a}sz extension for set functions. We also show that previous analyses of mixed-integer sets, such as the submodular knapsack polytope, can be extended to the approximate submodularity setting. Our work demonstrates that one may generalize many of the analytical tools used in submodular optimization into the approximate submodularity context. 

\end{abstract}

\keywords{Approximate submodularity, valid inequalities, set function extensions}

\section{Introduction} 
\label{section:Intro}
\vspace{-0.1in}

Exploiting structural properties in discrete optimization problems can lead to successful algorithms and heuristics. A classical property that is frequently used in discrete optimization is submodularity. Let $\Omega$ be a finite set of elements, and let $2^{\Omega}$ denote the power set of $\Omega$. A set function $f:2^{\Omega} \to \mathbb{R}$ is \textit{submodular} if for any $\cA \subseteq \cB \subseteq \Omega$ and $s \in \Omega \backslash \cB$, $f(\cB \cup \{s\}) - f(\cB) \leq f(\cA \cup \{s\}) - f(\cA)$. For some problems, submodularity provides guarantees for solution approaches such as the greedy algorithm. Recently, researchers have expanded algorithm analysis to \textit{approximately submodular} functions (e.g., \cite{Das2011,Horel2016,Zhou2016}). However, much of the initial focus on approximate submodularity has remained within performance guarantees for algorithms. In this paper, we propose approximate submodularity metrics to study multiple implications of approximate submodularity in discrete optimization, including the derivation of valid inequalities and properties of extensions on the unit hypercube. 
Our work applies to any nonnegative and monotonic set function, and our analyses often follow arguments similar to those of analogous results in the submodular context.

Continuous relaxations of problems are often used as direct approximation techniques for discrete optimization problems (e.g., solving the linear programming relaxation of a mixed-integer program) because they are easier to solve; these relaxations often have a polynomial-time algorithm. Extensions of set functions can transform discrete optimization problems into continuous optimization problems, for which efficient algorithms or approximation schemes may exist. An \textit{extension} of a set function $f:2^{\Omega} \to \mathbb{R}$ is a function $F:D \to \mathbb{R}$ such that $D \supset \mathbb{B}^{|\Omega|}$ and $F(x(\cS)) = f(\cS)$, for all $\cS \subseteq \Omega$, where $x(\cS)$ denotes the characteristic vector of the set $\cS$. We focus on extensions defined on the unit hypercube $[0,1]^{|\Omega|}$. Notably, the \emph{Lov\'{a}sz extension} \citep{Lovasz1983} for set functions is convex if and only if the set function is submodular; in this case, the Lov\'{a}sz extension is equal to the convex closure. The convex closure is difficult to compute in general; in contrast, computing the Lov\'{a}sz extension is comparatively simple, which makes submodularity a valuable property when considering solution methods that use the convex closure. \cutOut{Another useful extension of set functions is the multilinear extension. In particular, when the set function is submodular, its multilinear extension is up-concave (\thref{def:upConcave}), which has prompted researchers to use it in solution approaches for maximization problems.} We provide a new characterization that relates the approximate submodularity of a function with the approximate convexity of its Lov\'{a}sz extension.

Valid inequalities are crucial for solving mixed-integer programs as they can cut off solutions to relaxations so that the new problem's feasible region more closely approximates the convex hull \cite{Cornuejols2008}. The knapsack problem is one of the foundational problems in discrete optimization where researchers have studied its facial structure and valid inequalities (e.g., \cite{Atamturk2009,Balas1975}). In particular, Atamt\"{u}rk and Narayanan (2009) \cite{Atamturk2009} study valid inequalities for the submodular knapsack polytope, in which the constraint function is submodular. Submodular functions feature in the constraints of other optimization problems as well. Researchers have also studied mixed-integer programs with conic-quadratic constraints and objective functions where valid inequalities are derived by leveraging the submodularity of the objective and constraint functions \cite{Atamturk2008,Atamturk2020,Gomez2020}. Valid inequalities and outer approximations of the epigraphs of submodular and general set functions have also been studied \cite{AtamturkShort2020}. Our study is the first to use {\it{approximate}} submodularity to derive valid inequalities for mixed-integer sets defined by approximately submodular functions. 

Previous studies on approximate submodularity focus on performance bounds for greedy algorithms and other iterative selection approaches. Performance bounds have been produced using different notions of approximate submodularity where metrics with different properties can produce different bounds; trade-offs between additive and multiplicative bounds for the greedy algorithm performance on non-submodular functions are studied in \cite{Horel2016}, and \cite{Das2011} and \cite{Zhou2016} also define metrics that they use to propose greedy algorithm performance bounds for approximately submodular functions. However, the notion of approximate submodularity can also be used for generalizing results in other areas of discrete optimization, such as analyzing properties of continuous extensions and deriving valid inequalities, which can provide new insights for efficient solution methods. New methodological applications are still emerging, even outside of the greedy algorithm, in which approximate submodularity can extend the existing results that depend on submodularity. Our contributions are as follows:

\begin{itemize}
\item In Section \ref{sec:Metrics}, we study fundamental properties about our approximate submodularity metrics, which we use to show which operations preserve approximate submodularity. 
\item In Section \ref{sec:Extensions}, we derive results on the approximate convexity of the Lova\'{s}z extension of approximately submodular functions. 
\item We study mixed-integer sets defined by approximately submodular functions in Section \ref{sec:ValidInequalities}. We use the proposed metrics to adapt analogous analyses for the submodular setting, thus deriving new valid inequalities for cases when the set function is approximately submodular. 
\end{itemize}
We note that there are several cases in which our proofs are similar to those of analogous results in the submodular setting. Our primary message is that in a broad set of areas in discrete optimization, one can use approximate submodularity to generalize both classical and more recent results.

\vspace{-0.2in}
\section{Approximate Submodularity Metrics}\label{sec:Metrics}
\vspace{-0.1in}

In this section, we discuss various approximate submodularity metrics; the term ``metrics" is used loosely, as some are not subadditive and none are positive definite, both of which are part of the formal definition of a metric. However, the approximate submodularity metrics we discuss indicate a notion of distance to submodularity. In this work, if $\metric$ is an approximate submodularity metric, $\metric[f]$ is the metric value for $f$.

\vspace{-0.1in}
\subsection{Proposed Notions of Approximate Submodularity} \label{sec:ourMetrics}

We begin with the most general (global) metric, which is inspired directly from an equivalent definition of submodularity: $f(\cA \cup \cB) + f(\cA \cap \cB) \leq f(\cA) + f(\cB)$, for all $\cA, \cB \subseteq \Omega$.
\begin{definition}\thlabel{def:globalDistance}
Let $f: 2^{\Omega} \to \mathbb{R}$. Define the \emph{global submodularity distance} $\glob$ by $\glob[f] \coloneqq \max\limits_{\cA, \cB \subseteq \Omega} f(\cA \cup \cB) + f(\cA \cap \cB) - f(\cA) - f(\cB).$
\end{definition}
The global submodularity distance is a general purpose metric; we demonstrate its value in identifying operations that preserve approximate submodularity and proving general results about set functions (Section \ref{sec:PreserveApprox}). The remaining metrics are inspired by a characterization of 
increasing, submodular functions. 
\begin{lemma}\thlabel{def:Submod} \emph{(Edmonds 1970 \citep{Edmonds1970})}
Let $f: 2^{\Omega} \to \mathbb{R}$. Then $\myf$ is increasing and submodular if and only if for any $\cA, \cB \subseteq \Omega, s \in \Omega$, $f(\cA \cup \cB \cup \{s\}) - f(\cA \cup \cB) \leq f(\cA \cup \{s\}) - f(\cA)$.
\end{lemma}
Using \thref{def:Submod}, we present two metrics for approximate submodularity. 
\begin{definition}\thlabel{def:pairwise_viol}
Let $f: 2^{\Omega} \to \mathbb{R}$. Consider $\ell \in \{0,...,|\Omega| - 1\}, k \in \{0,...,|\Omega|\}$. The $(\ell,k)$\emph{-pairwise violation} of $\myf$ is defined as
$d^{\ell, k}[f] \coloneqq \max\limits_{\substack{\cA, \cB \subseteq \Omega, s \in \Omega \\|\cA| = \ell, |\cB| = k}} f(\cA \cup \cB \cup \{s\}) - f(\cA \cup \cB) - f(\cA\cup\{s\}) + f(\cA).$
\end{definition}
\noindent Thus, the pairwise violation represents the worst-case violation of the condition in \thref{def:Submod} given $\cA$ and $\cB$ with fixed cardinalities. In the context of a sensor placement problem, the $(\ell,k)$-pairwise violation captures the case in which a single sensor added to a sparse sensor network (given by $\cA$) creates a smaller marginal increase in information than when the same sensor is added to a denser network ($\cA \cup \cB$). Note that if $\myf$ is submodular, $d^{\ell, k}[f] \leq 0$ for all $\ell$ and $k$, and the reverse implication holds with the added condition of $\myf$ being monotonic increasing. Also, for any $f$, $d^{\ell,k}[f] \leq \glob[f]$.
 
\begin{definition}\thlabel{def:marginalViolation}
Let $f: 2^{\Omega} \to \mathbb{R}$. The \emph{marginal violation} of $f$ is defined as $\dknap[f] \coloneqq \max\big\{d^{\ell,k}[f] \setbar \ell \in \{0,\dots,|\Omega| - 1\}, k \in \{0,\dots,|\Omega|\}\big\}$.
\end{definition}
Note that $\dknap[f]$, also used in \cite{Krause10}, does not depend on set sizes. Although $\glob[f]$ and $\dknap[f]$ may be difficult to compute exactly in general, Section \ref{sec:PreserveApprox} details some operations that preserve approximate submodularity, with respect to $\glob[f]$, which 
enables one to bound $\glob[f]$ and $\dknap[f]$. In addition, in Appendix \ref{sec:facLocExample}, we present a generalized version of the uncapacitated facility location problem in which the objective function is approximately submodular and $\dknap[f]$ can be bounded analytically.

\vspace{-0.1in}
\subsection{Preserving Approximate Submodularity}\label{sec:PreserveApprox}
We prove some properties of our proposed approximate submodularity metrics (\thref{ApproxSubmodSubadditive}), as well as operations from which bounds or exact values of approximate submodularity metrics can be inferred immediately (\thref{submodPreserve}). The former compares properties of our approximate submodularity metrics to true metrics. The latter concept can be thought of as ``approximate submodularity preservation." Some of these results have analogs for submodular functions (for reference, see \cite{Bach2013}, \cite{Narayanan1997}, and \cite{NemhauserWolsey}), but others are specific to approximate submodularity. We let $\mathcal{F}$ (resp., $\mathcal{F}_{+}$) be the set functions (resp., that are nonnegative and increasing) over ground set $\Omega$.

\begin{theorem}\thlabel{ApproxSubmodSubadditive}
Consider a nonnegative, increasing set function $f: 2^{\Omega} \to \mathbb{R}$ and a metric of approximate submodularity $\metric: \mathcal{F}_{+} \to \mathbb{R}$ where $\metric$ is defined by any of the following: 
\begin{enumerate*}[label=(\Roman*)]
\item $\metric[f] = \glob[f]$,
\item $\metric[f] = \dknap[f]$, or
\item $\metric[f] = d^{\ell,k}[f], \text{ for some} \ \ell \in \{0,\dots,|\Omega| - 1\}, k \in \{0,\dots,|\Omega|\}$.
\end{enumerate*}
Then we have:
\begin{enumerate}[label=(\roman*)]
\item The function $\metric$ is sublinear. That is, $\metric$ is subadditive $($i.e., $\metric[f_{1}] + \metric[f_{2}] \geq \metric[f_{1} + f_{2}])$ and positively homogeneous with degree 1 $($i.e., $\alpha \metric[f] = \metric[\alpha f]$, for $\alpha \in \mathbb{R}_{+})$. \label{subaddClaim1}
\item If $f$ is not submodular, then for any $\epsilon \in [0, \metric[f])$, there does not exist a nonnegative, increasing, submodular function $g: 2^{\Omega} \to \mathbb{R}$ such that $||g - f||_{\infty} < \frac{\epsilon}{4}$. \label{subaddClaim2}
\end{enumerate}
\end{theorem}

The contrapositive of Claim \ref{subaddClaim2} of \thref{ApproxSubmodSubadditive} can be read as a necessary condition, which can, in some cases, remove the need for testing whether any function near $f$ is submodular (e.g., \cite{Seshadri}). Although our notions of approximate submodularity are not ``metrics" in the analytical sense, \thref{ApproxSubmodSubadditive} proves that they are sublinear. Sublinear functions are well studied in the literature and are the ``next simplest convex functions" after affine functions \citep{Urruty2001}. All metrics are sublinear. Subadditivity and positive homogeneity independently have multiple implications. They can be used to verify that a function $f_{1} + f_{2}$ (or $\alpha f$, for $\alpha \in \mathbb{R}_{+}$) satisfies conditions in hypotheses of results in Sections \ref{sec:Extensions}--\ref{sec:ValidInequalities}. We remark that subadditivity (and hence, sublinearity) is not a trivial property of approximate submodularity metrics in the literature; e.g., submodularity ratio proposed in Das and Kempe (2011) \cite{Das2011} is not subadditive.

We can relate our metrics to asymmetric seminorms, which share more properties with analytical metrics. 
\begin{definition}\thlabel{def:asymSemiNorm} \emph{(Cobza\c{s} 2013 \citep{Cobzas2012})}
A function $\metric:\mathbb{R}^{n} \to \mathbb{R}$ is an asymmetric seminorm if it is nonnegative, positively homogeneous, and subadditive.
\end{definition}
\begin{corollary}\thlabel{induceASN}
Define $\glob_{+}:\mathcal{F} \to \mathbb{R}$ by $\glob_{+}[f] \coloneqq \max\{0, \glob[f]\}$. Then $\glob_{+}$ is an asymmetric seminorm on $\mathcal{F}$. 
\end{corollary}

We provide some examples in which functions induced by an approximately submodular function inherit approximate submodularity. Denote the complement of $\cS \subseteq \Omega$ by $\cS^{c}$. Given a normalized set function $f$, define $f_{1}, f_{2}: 2^{\Omega} \to \mathbb{R}$ by $f_{1}(\cS) = f(\cS^{c}), f_{2}(\cS) = f(\cS) + f(\cS^{c}) - f(\Omega)$; thus, $f_{2}$ is a symmetric, nonnegative function. Given $\cA \subseteq \Omega,$ define $f_{\cA}:2^{\Omega \backslash \cA}\to\mathbb{R}$ by $f_{\cA}(\cS) = f(\cA \cup \cS)$. Given a factor $q$ of $|\Omega|$, let $\Omega(q) = \{1,\dots,\frac{|\Omega|}{q}\}$, $\cS(i) = \{(i -1)q + 1,\dots,iq\}$, for all $i \in \Omega(q)$, and $f_{q}:2^{\Omega(q)} \to \mathbb{R}$ be defined by $f_{q}(\cS) = f(\bigcup\limits_{i \in \cS} \cS(i))$. Finally, let $g:2^{\Omega} \to \mathbb{R}$ be a \textit{modular} function ($g$ and $-g$ are submodular), and define the convolution of $f$ and $g$ as $f\circledast g(\cS) = \min\limits_{\cZ \subseteq \cS} f(\cZ) + g(\cS \backslash \cZ)$. Note $f\circledast g(\cS) = g \circledast f(\cS)$. 

\begin{proposition}\thlabel{submodPreserve}
Given $f:2^{\Omega} \to \mathbb{R}$ and the corresponding functions $f_{1}, f_{2}, f_{\cA},$ and $f_{q}$,
we have:
\begin{enumerate*}[label=(\roman*)]
\item $\glob[f] = \glob[f_{1}]$. \label{submodPreserve1}
\item $2\glob[f] \geq \glob[f_{2}]$. \label{submodPreserve2}
\item $\glob[f] \geq \glob[f_{\cA}]$ \label{submodPreserve3}
\item $\glob[f] \geq \glob[f_{q}]$. \label{submodPreserve4}
\item $\glob[f] \geq \glob[f\circledast g]$. \label{submodPreserve5}
\end{enumerate*}
\end{proposition}
\thref{submodPreserve} can be used in a fashion similar to \thref{ApproxSubmodSubadditive}. Also, \thref{submodPreserve}\ref{submodPreserve4} can provide guarantees on a greedy algorithm that selects among prescribed subsets of elements. Our proof for \thref{submodPreserve}\ref{submodPreserve5} follows similar arguments to that of the submodular case \citep{Narayanan1997}. Note that we slightly abuse notation in \thref{submodPreserve}\ref{submodPreserve3}--\ref{submodPreserve4} as the domains of $f_{\cA}$ and $f_{q}$ are not $2^{\Omega}$. 

\section{Extensions of Approximately Submodular Functions}\label{sec:Extensions}
\vspace{-0.1in}

Next, we study extensions of approximately submodular functions. For the remainder of this paper, we consider functions that are monotone increasing and normalized ($f(\emptyset) = 0$) and focus on analysis based on the marginal violation $D$ because we consider problems of this form, namely the approximately submodular knapsack and packing problems and the generalized uncapaciated facility location problem, in Section~\ref{sec:ValidInequalities} and Appendix~\ref{sec:facLocExample}, respectively. Other examples for monotone submodular optimization can be found in \cite{Krause2014}, which may have relevant approximately submodular analogs. Given a set function $f:2^{\Omega} \to \mathbb{R}$, an \textit{extension} of $f$ over $[0,1]^{|\Omega|}$ is a function $F:[0,1]^{|\Omega|} \to \mathbb{R}$ such that $F(x(\cS)) = f(\cS)$, for all $\cS \subseteq \Omega$, where $x(\cS)$ is the characteristic vector of $\cS$. Our main result in this section is that the Lov\'{a}sz extension is approximately convex (\thref{LovaszApprox}) when the underlying set function is approximately submodular. A main component of multiple key results in this section is the marginal violation $\dknap$ (\thref{def:marginalViolation}). Other works on extensions of set functions include \cite{Iyer2015}, \cite{Lovasz1983}, and \cite{Murota1998}.

The Lov\'{a}sz extension of a set function $f: 2^{\Omega} \to \mathbb{R}$ is defined by $F^{L}: [0,1]^{|\Omega|} \to \mathbb{R}$ such that $F^{L}(x) \coloneqq \sum\limits_{k = 0}^{|\Omega|} \lambda_{k}f(\cC_{k})$, where $\emptyset = \cC_{0} \subset \cC_{1} \subset \cdots \subset \cC_{|\Omega|} = \Omega$ is a chain such that $\sum\limits_{k = 0}^{|\Omega|}\lambda_{k}x(\cC_{k}) = x$, with $\sum\limits_{k = 0}^{|\Omega|}\lambda_{k} = 1$, and $\lambda \geq 0$. It is well known that by defining a permutation $(\pi_{1},\dots,\pi_{|\Omega|})$ such that $x_{\pi_{1}} \geq x_{\pi_{2}}\cdots\geq x_{\pi_{|\Omega|}}$, $\cC_{0} = \emptyset, \cC_{k} = \cC_{k-1} \cup \{\pi_{k}\}$, for $k \in \{1,\dots,|\Omega|\}$, the Lov\'{a}sz extension is equivalently defined as $F^{L}(x) = \sum\limits_{k = 1}^{|\Omega|} x_{\pi_{k}}(f(x(\cC_{k})) - f(x(\cC_{k-1})))$ (see \cite{Bach2013}). The convex closure of $f$ is the unique convex function $F^{C}:[0,1]^{|\Omega|}\to\mathbb{R}$ such that $F^{C}(x(\cS)) \leq f(\cS)$ for all $\cS \subseteq \Omega$ and $F^{C}(x) \geq G(x)$ for any other convex understimator $G:[0,1]^{|\Omega|}\to\mathbb{R}$ of $f$. Lov\'{a}sz (1983) \cite{Lovasz1983} shows that $f$ is submodular if and only if $F^{L}$ is convex; in fact, in this special case, the convex closure and the Lov\'{a}sz extension are equal ($F^{C} = F^{L}$). This property is useful in that the convex closure is generally difficult to compute in comparison to the Lov\'{a}sz extension. Although the Lov\'{a}sz extension does not equal the convex closure when $f$ is not submodular, we prove a generalized result when $f$ is approximately submodular. We remark that Halabi and Jegelka (2019) \cite{Halabi2019} also study the Lov\'{a}sz extension of non-submodular functions, including its subgradients, in the context of convex optimization solution approaches.

Consider the following linear program parametrized by $x \in [0,1]^{|\Omega|}$:
\begin{equation}\label{closureLPDual}
\begin{aligned}
V(x) = \min\limits_{y}\left\{ \sum\limits_{\cS \subseteq \Omega}f(\cS)y(\cS) \ \Bigg\vert \ \sum\limits_{\cS\ni s} y(\cS) = x_{s}, \forAll s \in \Omega,
\sum\limits_{\cS \subseteq \Omega} y(\cS) = 1,
y \geq 0\right\}.
\end{aligned}
\end{equation}

\begin{proposition}{\emph{(Bach 2013 \citep{Bach2013})}}\thlabel{convexClosureLP}
For $f:2^{\Omega}\to\mathbb{R}$ with $f(\emptyset) = 0$, we have $V(x) = F^{C}(x)$, for all $x \in [0,1]^{|\Omega|}$.
\end{proposition}

Given a permutation $(\pi_{1},\cdots,\pi_{|\Omega|})$, define $\cS_{0}^{\pi} \coloneqq \emptyset \subset \cS_{1}^{\pi} \coloneqq \{\pi_{1}\}\cdots\cS_{k}^{\pi} \coloneqq \{\pi_{1},\dots,\pi_{k}\} \cdots\subset \cS_{|\Omega|}^{\pi} \coloneqq \Omega$. Define the set
$\Gamma(f) \coloneqq \{\gamma \in \mathbb{R}^{|\Omega|} \setbar \exists \text{ permutation } \pi \text{ such that } \gamma_{\pi_{i}} = f(\cS^{\pi}_{i}) - f(\cS^{\pi}_{i-1}), \forAll i \in \Omega\}.$

\begin{definition}
A function $F: [0,1]^{|\Omega|} \to \mathbb{R}$ is \emph{$\epsilon$-approximately convex}, if $F(\lambda x + (1 - \lambda)y) \leq \epsilon + \lambda F(x) + (1 - \lambda)F(y)$, for any $\lambda \in [0,1]$.
\end{definition}

\begin{theorem}\thlabel{LovaszApprox}
For any increasing set function $f:2^{\Omega} \to \mathbb{R}$ such that $f(\emptyset) = 0$, 
\begin{align*}
F^{L}(x) &\leq \max\limits_{\gamma \in \Gamma(f)} \sum\limits_{s \in \Omega} \gamma_{s}x_{s}
\leq F^{C}(x) + |\Omega|\dknap[f]
\leq F^{L}(x) + |\Omega|\dknap[f] 
\leq \max\limits_{\gamma \in \Gamma(f)} \sum\limits_{s \in \Omega} \gamma_{s}x_{s} + |\Omega|\dknap[f].
\end{align*}
Hence, $F^{L}(x) \geq F^{C}(x) \geq F^{L}(x) - |\Omega|\dknap[f]$, and $||F^{L} - F^{C}||_{\infty} \leq |\Omega|\dknap[f]$. Moreover, $F^{L}$ is $|\Omega|\dknap[f]$-approximately convex.
In addition, if for some $\epsilon > 0$, $F^{L}$ is $\epsilon$-approximately convex,  then $\dknap[f] \leq \epsilon$.
\end{theorem}

\thref{LovaszApprox} states that the approximate submodularity of $f$ implies the approximate convexity of $F^{L}$ and vice-versa. The proof of \thref{LovaszApprox} uses the well-known linear program \eqref{closureLPDual}, but a key difference is that we construct feasible primal-dual solutions with a duality gap due to the generalization to approximate submodularity.

Next, we consider the case in which there exists a submodular function $g$ close to $f$. In this case, we show that the Lov\'{a}sz extension of $g$ approximates the Lov\'{a}sz extension of $f$.
\begin{proposition}\thlabel{SubmodToConvexLovasz}
Given set functions $f,g:2^{\Omega} \to \mathbb{R}$, where $f(\emptyset) = g(\emptyset) = 0$, and their respective Lov\'{a}sz extensions $F^{L},G^{L}: [0,1]^{|\Omega|} \to \mathbb{R}$, $||F^{L} - G^{L}||_{\infty} = ||f - g||_{\infty}$.
\end{proposition}

Thus, the approximating function $g$ (which may be submodular) can lead to approximation methods in the discrete domain or over the hypercube using convex optimization methods.

\section{Valid Inequalities of Polyhedra Associated With Approximately Submodular Functions}\label{sec:ValidInequalities}
\vspace{-0.1in}

We use approximate submodularity metrics from Section \ref{sec:ourMetrics} to derive valid inequalities for some mixed-integer sets. Our analyses are similar to analogs in submodular analysis \cite{Atamturk2008,Atamturk2009,AtamturkShort2020} with additional details to generalize to approximate submodularity. 

\subsection{Epigraph Inequalities}\label{sec:epigraph}
First, we study the epigraphs of set functions. These mixed-integer sets can be useful when minimizing a submodular function \citep{Atamturk2020}. We consider the case when the function is approximately submodular. Let $\phi: \mathbb{R} \to \mathbb{R}_{+}$ be increasing, and for any $\tau \in \mathbb{R}_{+}$, let $F_{\tau}:[0,1]^{|\Omega|} \to \mathbb{R}_{+}$ be defined by $F_{\tau}(x) = \phi(\tau + \sum\limits_{i \in \Omega}c_{i}x_{i})$, where $c \in \mathbb{R}^{|\Omega|}_{+}$. Thus, $F_{\tau}$ is increasing. Consider the mixed-integer feasible region
$H_{\mathbb{B}} = \left\{(x,z) \in \mathbb{B}^{|\Omega|} \times \mathbb{R}_{+} \setbar F_{\sigma}(x) \leq z\right\},$
where $\sigma \geq 0,$ and $c \in \mathbb{R}^{|\Omega|}_{+}$. 
Define the set function $f_{\tau}:2^{\Omega} \to \mathbb{R}$ by $f_{\tau}(\cS) = F_{\tau}(x(\cS))$. Notice that $f_{\tau}$ is increasing and $f_{\tau}(\emptyset) = 0$ if and only if $\phi(\tau) = 0$. Therefore, define $g_{\tau}:2^{\Omega} \to \mathbb{R}$ by $g_{\tau}(\cS) = f_{\tau}(\cS) - \phi(\tau)$, which is normalized, $g_{\tau}(\emptyset) = 0$, and is increasing; hence, it is also nonnegative. Note that $\dknap[f_{\tau}] = \dknap[g_{\tau}]$. We denote the Lov\'{a}sz extension of $g_{\tau}$ by $G_{\tau}^{L}$.

For any $\gamma \in \Gamma(g_{\tau})$---i.e., $\gamma_{\pi_{k}} = f_{\tau}(\cS^{\pi}_{k}) - f_{\tau}(\cS^{\pi}_{k-1}) = \phi(\tau + \sum\limits_{i = 1}^{k}c_{\pi_{i}}) - \phi(\tau + \sum\limits_{i = 1}^{k}c_{\pi_{i-1}})$ for some permutation $\pi$ of $(1,\dots,|\Omega|)$---consider the following inequality:
\begin{align}
\sum\limits_{s \in \Omega}\gamma_{s}x_{s} \leq z - \phi(\tau). \label{ASP_Inequality}
\end{align}

When $\phi$ is the square root function, then $f_{\tau}$ is a submodular set function, and Atamt\"{u}rk and Narayanan (2008) \cite{Atamturk2008} show that inequality \eqref{ASP_Inequality} is valid for the convex hull $\conv(H_{\mathbb{B}})$ for $\tau = \sigma$. In fact, along with the variable bounds, such inequalities describe $\conv(H_{\mathbb{B}})$. In the more general case, where $\phi$ is such that $f_{\tau}$ is approximately submodular, we show that similar inequalities are still valid for $\conv(H_{\mathbb{B}})$. 

\begin{lemma}\thlabel{generalInequality}
For any $\gamma \in \Gamma(g_{\tau}), \cS \subseteq \Omega,$ we have $-|\Omega|\dknap[f_{\tau}] + \sum\limits_{s \in \Omega} \gamma_{s}(x(\cS))_{s} \leq f_{\tau}(\cS) - \phi(\tau)$.
\end{lemma}

\begin{proposition}\thlabel{validEpigraphBinary}
For any $\gamma \in \Gamma(g_{\sigma})$, the following inequality is valid for $\conv(H_{\mathbb{B}})$:
\begin{align}\label{validSigmaInequality}
-|\Omega|\dknap[f_{\sigma}] + \sum\limits_{s \in \Omega}\gamma_{s}x_{s} \leq z - \phi(\sigma).
\end{align} 
\end{proposition}

\thref{validEpigraphBinary} illustrates what is lost between submodularity and approximate submodularity in deriving valid inequalities in this setting. When $f_{\sigma}$ and $g_{\sigma}$ are approximately submodular, $\dknap[f_{\tau}] > 0$ may lead to looser valid inequalities. Our proof of \thref{validEpigraphBinary} follows arguments similar to those of Atamt\"{u}rk and Narayanan (2008) \cite{Atamturk2008}, who establish the result when $\phi$ is the square root function.

Next, we consider the epigraph of a general, increasing, nonnegative, approximately submodular function $f$, $H_{f} \coloneqq \mathrm{conv}(\{(x,z) \in \mathbb{R}^{|\Omega|}\times\mathbb{R} \setbar f(x(\cS)) \leq z\})$. Consider the associated polyhedron $P_{f} \coloneqq \{\gamma \in \mathbb{R}^{|\Omega|} \setbar \sum\limits_{s \in \Omega} \gamma_{s} \leq f(\cS), \forAll \cS \subseteq \Omega\}$. We refer to the variable bounds as trivial inequalities of $H_{f}$.

\begin{proposition}\emph{(Atamt\"{u}rk and Narayanan 2020 \citep{AtamturkShort2020})}\thlabel{atamturkOuter}
\begin{enumerate}
\item Any nontrivial facet-defining inequality $\sum\limits_{s \in \Omega} \gamma_{s}x_{s} \leq \alpha z + \gamma_{0}$ for $H_{f}$ satisfies $\gamma_{0} \geq 0$ and $\alpha = 1$ (up to scaling).
\item The inequality $\sum\limits_{s \in \Omega}\gamma_{s}x_{s} \leq z$ is valid for $H_{f}$ if and only if $\gamma \in P_{f}$.
\item The inequality $\sum\limits_{s \in \Omega}\gamma_{s}x_{s} \leq z$ is facet-defining for $H_{f}$ if and only if $\gamma$ is an extreme point of $P_{f}$.
\end{enumerate}
\end{proposition}
Atamt\"{u}rk and Narayanan (2020) \cite{AtamturkShort2020} prove that nontrivial facets of $H_{f}$ are homogeneous. We establish a similar result for approximately submodular functions.

\begin{proposition}\thlabel{polarIneq}
Let $f:2^{\Omega} \to \mathbb{R}$ be increasing with $f(\emptyset) = 0$. Suppose $\gamma \in \mathbb{R}^{|\Omega|}$ and
\begin{align}
\sum\limits_{s \in \Omega}\gamma_{s}x_{s} \leq z + |\Omega|\dknap[f]+\gamma_{0} \label{ineqGeneralOut}
\end{align}
defines a nontrivial facet of $H_{f}$. Let $\bar{f}:2^{\Omega} \to \mathbb{R}$ be defined by $\bar{f}(\emptyset) = 0, \bar{f}(\cS) = f(\cS) + |\Omega|\dknap[f] + \gamma_{0}$, for all nonempty $\cS \subseteq \Omega$, and suppose $\gamma \in \Gamma(\bar{f})$. Then, $\gamma_{0} \leq 0$.
\end{proposition}

The proof of \thref{polarIneq} proceeds similarly to that of the submodular case in \cite{AtamturkShort2020}, with some additional steps to account for the approximate submodularity generalization. This includes bounding $\max\limits_{\gamma \in \Gamma[f]} \sum\limits_{s \in \cS} \gamma_{s}$ using the marginal violation $\dknap$. Given the conditions in the hypothesis of \thref{polarIneq}, the constant term $|\Omega|\dknap[f] + \gamma_{0}$ is bounded below by 0 and above by $|\Omega|\dknap[f]$; when $f$ is submodular, the condition $\gamma \in \Gamma(\bar{f})$ is implied, $\dknap[\bar{f}] = 0$, and the nontrivial facets are homogeneous. We also remark that Atamt\"{u}rk and Narayanan (2020) \cite{AtamturkShort2020} provide valid inequalties for general set functions, but these rely on a submodular-supermodular decomposition of $f$.

\subsection{Knapsack Inequalities} \label{sec:knapIneq}
Consider the polytope $X = \text{conv}\{x \in \mathbb{B}^{n} \setbar f(\cS(x)) \leq b\}$, where $\cS(x)$ is the subset of $\Omega$ characterized by the binary vector $x$ and $f$ is a set function. When $f$ is submodular, nonnegative, and increasing, $X$ is known as the submodular knapsack polytope \citep{Atamturk2009}; a special case is the well-known linear knapsack set, and optimizing over it is NP-hard \citep{Karp}. We consider the case where $f$ is approximately submodular, nonnegative and increasing. Thus, we call the set $X$ an \textit{approximately submodular knapsack set}. Our focus in this subsection is on deriving valid inequalities for this set. Some facets for 0-1 polytopes established by Atamt\"{u}rk and Narayanan (2009) \cite{Atamturk2009} apply in our setting; we list them in \thref{thm:easyInequalities} in the appendix. 

\begin{definition}\thlabel{def:Covers}
\text{}
\begin{enumerate}
\item The subset $\cS \subseteq \Omega$ is a \emph{cover} for $X$ if $f(\cS) > b$ and is \emph{minimal} if $f(\cS \backslash \{s\}) \leq b$ for all $s \in \cS$.
\item Let $\pi = (\pi_{1},\dots,\pi_{|\Omega\backslash \cS|})$ be a permutation of $\Omega \backslash \cS$. Let $U_{\pi}(\cS) = \{\pi_{j} \in \Omega \backslash \cS \setbar f(\cS \cup \{\pi_{1},\dots,\pi_{j}) - f(\cS\cup\{\pi_{1},\dots,\pi_{j-1}\}) \geq f(\{s\}), \forAll s \in \cS\}$.
The \emph{set-extension} of $\cS \subseteq \Omega$ with respect to $\pi$ is denoted by $E_{\pi}(\cS) = \cS \cup U_{\pi}(\cS)$. 
\end{enumerate}
\end{definition}

\thref{Atamturk09ExtendedV2} extends the result of Proposition 5 in Atamt\"{u}rk and Narayanan (2009) \cite{Atamturk2009} for submodular knapsack problems into the approximately submodular context.

\begin{proposition}\thlabel{Atamturk09ExtendedV2}
If $\cS \subseteq \Omega$ is a cover for $X$, the extended cover inequality $\sum\limits_{s \in E_{\pi}(\cS)} x_{s} \leq |\cS| - 1$ is valid for $X$ if $f(\cS) > (|\cS| + |U_{\pi}(\cS)|)\dknap[f] + b$. In addition, the inequality defines a facet of $\{x \in X \setbar x_{s} = 0, \forAll s \not\in E_{\pi}(\cS)\}$ if $\cS$ is also a minimal cover and for each $s \in U_{\pi}(\cS)$, there exist $t_{s}, u_{s} \in \cS$ such that $t_{s} \neq u_{s}$, and $f(\cS \cup \{s\} \backslash \{t_{s},u_{s}\}) \leq b$.
\end{proposition}

We observe from \thref{Atamturk09ExtendedV2} that in adapting the result for approximate submodularity, we add a condition for the extended cover inequality to be valid. The proof follows similar steps as those in Atamt\"{u}rk and Narayanan (2009) \cite{Atamturk2009}, except it accounts for violated submodularity inequalities.

\subsection{Illustrative Example} \label{sec:numerical}
In this section we provide an example in which we derive valid inequalities for the approximately submodular knapsack polytope. We apply the derived valid inequalities (\thref{Atamturk09ExtendedV2}) and show that they can be used as a tool in the process of finding integer solutions when optimizing a linear function over the polytope.

Let $\Omega$ be a set of elements, $u, w \in \mathbb{R}^{|\Omega|}_{+} \backslash \{0\}$, $p > 1$. Let $G:[0,1]^{|\Omega|} \to \mathbb{R}$ and $H:\mathbb{R} \to \mathbb{R}$, where $G(x) = w^{\top}x$ and $H(z) = z^{p}$ if $z \geq 0$ and $0$ otherwise. Define $F: [0,1]^{|\Omega|} \to \mathbb{R}$ and $f:2^{\Omega} \to \mathbb{R}$ by $F(x) = u^{\top}x + H(G(x))$ and $f(\cS) = F(x(\cS))$.
Because $p > 1$, $H$ is convex and increasing. Also, $G$ is a linear function, so $F$ is convex and $f$ is supermodular. Moreover, $f$ is increasing and nonnegative but not submodular. 

\begin{proposition}\thlabel{knapsackFunctionBound}
We have $\dknap[f] \leq p||w||_{1}^{p-1}||w||_{\infty}$. 
\end{proposition}

\thref{knapsackFunctionBound} gives a bound on $\dknap[f]$ that can remove the need to compute $\dknap[f]$ directly, which helps verify whether \thref{Atamturk09ExtendedV2} applies to inequalities of the form $f(\cS) > (|\cS| + |U_{\pi}(\cS)|)\dknap[f] + b,$ for some $\cS \subseteq \Omega$. 

We provide an example instance in which we optimize a linear function over the integer hull of the approximately submodular knapsack polytope. In particular, we show that by adding our valid inequalities, it is possible to obtain an integral solution from the continuous relaxation. Let $\Omega = \{1,2,\dots,6\}$, $u = [9,9,9,9,8.85,0]^{\top},$ $w = [0,0,0,0,1,1]^{\top},$ $p = 1.1$, and define $f, F, G,$ and $H$ as stated above. Also let $c = [3,3,3,3,2,2]^{\top}, b = 28.3$. Define the following instance of an approximately submodular knapsack problem \eqref{ASK} written as a binary program:
\begin{equation}\tag{ASK}\label{ASK}
\begin{aligned}
z_{\mathrm{ASK}} = \max\limits\left\{\sum\limits_{s \in \Omega} c_{s}x_{s} \setbar F(x) \leq b, x \in \mathbb{B}^{|\Omega|}\right\}.
\end{aligned}
\end{equation}

The continuous relaxation of \eqref{ASK} was solved using Gurobi 9.1.1 \citep{Gurobi} through the Gurobipy python interface (Python version 3.6.8) using a piecewise approximation of the nonlinear function with maximum absolute error of .001, with an optimal solution of $[0.0\overline{33}, 1, 1, 1, 0, 1]^{\top}$ and objective value $11.1$. Observe that $\cS = \{1,2,3,4\}$ is a (minimal) cover. Consider the permutation of $(5,6)$ $(\pi_{1} = 5, \pi_{2} = 6)$. Then $f(\cS \cup \{\pi_{1}\}) - f(\cS) = 9.85 > 9 = f(\{s\}),$ and $f(\cS \cup \{\pi_{1},\pi_{2}\}) - f(\cS \cup \{\pi_{1}\}) \approx 1.14 < f(\{s\})$, for all $s \in \cS$. Hence, $U_{\pi}(\cS) = \{5\}$ and $E_{\pi}(\cS) = \{1,2,\dots,5\}$. By \thref{knapsackFunctionBound}, $\dknap[f] \leq p||w||_{1}^{p-1}||w||_{\infty} \approx 1.18$, which implies $f(\cS) = 36 > 1.18|\Omega| + 28.3 = 7.08 + 28.3 \geq |\Omega|\dknap[f] + b \geq (|\cS| + |U_{\pi}(\cS)|)\dknap[f] + b;$ thus, \thref{Atamturk09ExtendedV2} implies $\sum\limits_{s = 1}^{5} x_{s} \leq 3$ is a valid inequality for \eqref{ASK}. Solving the relaxation of \eqref{ASK} with this valid inequality yields an optimal solution of $[1,1,1,0,0,1]^{\top}$ with an objective value of $11$. Thus, this solution is optimal for \eqref{ASK}. 

In general, $f(\cS) > (|\cS| + |U_{\pi}(\cS)|)\dknap[f] + b$ does not hold, so not every extended cover inequality is valid. 

\subsection{Randomly Generated Instances} \label{sec:randInstances}
We further explore the utility of the valid inequalities presented in \thref{Atamturk09ExtendedV2} by solving 40 randomly generated instances with and without a selection of these valid inequalities. The instances are larger versions of the example instance illustrated in Section \ref{sec:numerical}; each instance contains $q = 10$ approximately submodular knapsack constraints, and we refer to this problem as an \textit{approximately submodular packing problem}. Table \ref{tab:paramTable} shows the lower and upper bounds  of the uniform distributions of the parameters of the problem instances.

Because there were $n = 120$ decision variables, it was not practical to examine all possible valid inequalities for each instance. Instead, for each knapsack constraint, we generated 90 sets of $\cS$ where $|\cS|\in\{5,6,7,8\}$ and searched for valid inequalities based on these sets. In particular, for five repetitions, we randomly selected a subset $\cT \subset \cS^{c}$ where $|\cT|=8$ and determined if an extended cover valid inequality could be generated. To aid the search for valid inequalities near the best feasible solutions, the probability that variable $i$ was added to set $\cS$ was given by $\frac{c_{i}}{\sum_{j} c_{j}}$. All generated inequalities were added to the formulation before the start of the solve. We note that optimizing the inequality generation process is outside the scope of this study; thus, we do not include the generation time in our results and leave this subject to future research. All problems were solved using Gurobi 9.1.1 \citep{Gurobi} through the Gurobipy Python interface (Python version 3.6.8) on a Linux machine with a 20-core 2.4 GHz processor and 256 GB memory.

\begin{figure}
\centering
\includegraphics[width=.8\linewidth]{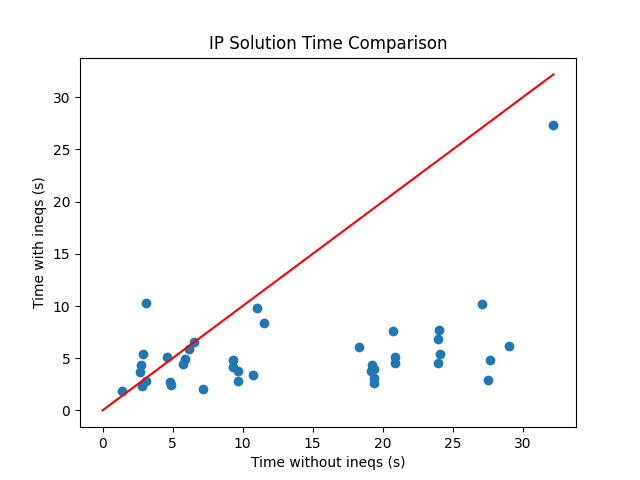}
\caption{Solution times (with and without inequalities) in seconds for 40 instances of approximately submodular packing problems. Points below the 
diagonal line (where the time ratio equals 1) indicate reduced solution times by adding extended cover inequalities.}
\label{fig:timeRatioScatter}
\end{figure}

Overall, we observe that the extended cover valid inequalities improve the solution time for most of the instances (see Figure \ref{fig:timeRatioScatter}). The cumulative time to solve all of the instances with the extended cover inequalities was less than half that without the additional inequalities. Table \ref{tab:InstanceSolveTime} in the appendix provides the solution times and the number of added inequalities for each instance.

\begin{table}
\centering
\begin{tabular}{|c|c|c|c|}
\hline
\bfseries Parameter & \bfseries LB & \bfseries UB & Set Value\\
\hline
$w_{i}$ & .2 & .6 & $\cdot$\\
\hline
$u_{i}$ & 1 & 11 & $\cdot$\\
\hline
$c_{i}$ & 8.5 & 11.5 & $\cdot$\\
\hline
$b$ & 42 & 48 & $\cdot$\\
\hline
$p$ & 1.05 & 1.2 & $\cdot$\\
\hline
$q$ & $\cdot$ & $\cdot$ & 10\\
\hline
$n$ & $\cdot$ & $\cdot$ & 120\\
\hline
\end{tabular}
\caption{Approximately submodular packing problem instance parameter values for $w, u, c, b,$ and $p$ were randomly generated from uniform distributions with the above lower and upper bounds. All problem instances had exactly 120 decision variables.}\label{tab:paramTable}
\end{table}

\section{Conclusion}
The value of submodularity in discrete optimization has long been established. Recently, notions of approximate submodularity have been applied to the greedy algorithm and similar approaches. We introduce new approximate submodularity metrics that have broad applicability in discrete optimization. 
We derive fundamental properties about our metrics, including which set function operations preserve approximate submodularity. We establish connections between our notions of approximate submodularity  and the approximate convexity of the Lov\'{a}sz extension.  \label{pgref:endConclusion}
Our approximate submodularity metrics can directly extend analyses in areas such as valid inequality derivations. Our numerical results show that valid inequalities derived based on the proposed metrics reduce solution times for approximately submodular packing problems in general. Optimizing the generation of valid inequalities for these problems and comparing their performances to prior works remains our future work.

\section*{\large Data availability statement}
The implementable instances used in this study are available at:\\ \url{https://bitbucket.org/tayoajayi/approxsubmodinstances2021/src/main/}.

\vspace{-0.2in}

\section*{\large Acknowledgements}
\vspace{-0.1in}
The authors would like to thank the anonymous referees and associate editor, Seth Brown, David Mildebrath, Logan Smith, and Silviya Valeva of Rice University for their helpful comments. This research was funded by National Science Foundation grants CMMI-1826297 and CMMI-1826323.

\vspace{-0.2in}

\bibliography{OPTL-D-21-00133R1_spring2022.bib}
\bibliographystyle{plain}

\newpage
\appendix

\section{Omitted Proofs}
\begin{repeattheorem}{ApproxSubmodSubadditive}
Consider a nonnegative, increasing set function $f: 2^{\Omega} \to \mathbb{R}$ and a metric of approximate submodularity $\metric: \mathcal{F}_{+} \to \mathbb{R}$ where $\metric$ is defined by any of the following: 
\begin{enumerate*}[label=(\Roman*)]
\item $\metric[f] = \glob[f]$,
\item $\metric[f] = \dknap[f]$, or
\item $\metric[f] = d^{\ell,k}[f], \text{ for some} \ \ell \in \{0,\dots,|\Omega| - 1\}, k \in \{0,\dots,|\Omega|\}$.
\end{enumerate*}
Then we have:
\begin{enumerate}[label=(\roman*)]
\item The function $\metric$ is sublinear. That is, $\metric$ is subadditive $($i.e., $\metric[f_{1}] + \metric[f_{2}] \geq \metric[f_{1} + f_{2}])$ and positively homogeneous with degree 1 $($i.e., $\alpha \metric[f] = \metric[\alpha f]$, for $\alpha \in \mathbb{R}_{+})$. \label{subaddClaim1}
\item If $f$ is not submodular, then for any $\epsilon \in [0, \metric[f])$, there does not exist a nonnegative, increasing, submodular function $g: 2^{\Omega} \to \mathbb{R}$ such that $||g - f||_{\infty} < \frac{\epsilon}{4}$. \label{subaddClaim2}
\end{enumerate}
\end{repeattheorem}

\myProof \label{pgref:caseIII}
For both claims, we prove case (III); proofs for other cases are similar. Let $f_{j}:2^{\Omega} \to \mathbb{R}_{+}$ be increasing set functions for $j \in \{1,2\}$.  

Consider $\metric(\cdot) = d^{\ell,k}(\cdot), \ell \in\{0,\dots,|\Omega| - 1\}, k \in \{0,\dots,|\Omega|\}.$ Observe that
\begin{align*}
\metric[f_{1}+f_{2}] &= \max\limits_{\substack{\cA, \cB \subseteq \Omega, s \in \Omega \\|\cA| = \ell, |\cB| = k}} \sum\limits_{j = 1}^{2} \Big(f_{j}(\cA \cup \cB \cup \{s\}) - f_{j}(\cA \cup \cB) - f_{j}(\cA\cup\{s\}) + f_{j}(\cA)\Big)\\
&\leq \sum\limits_{j = 1}^{2} \left(\max\limits_{\substack{\cA, \cB \subseteq \Omega, s \in \Omega \\|\cA| = \ell, |\cB| = k}} \Big(f_{j}(\cA \cup \cB \cup \{s\}) - f_{j}(\cA \cup \cB) - f_{j}(\cA\cup\{s\}) + f_{j}(\cA)\Big)\right)\\
&= \metric[f_{1}] + \metric[f_{2}],
\end{align*}
which proves subadditivity.

Let $(\cA^{*}, \cB^{*}, s^{*}) \in \arg\max\limits_{\substack{\cA, \cB \subseteq \Omega, s \in \Omega \\|\cA| = \ell, |\cB| = k}} \Big( f(\cA \cup \cB \cup \{s\}) - f(\cA \cup \cB) - f(\cA\cup\{s\}) + f(\cA) \Big)$. Then, for any $\alpha \in \mathbb{R}_{+}$, we have $(\cA^{*}, \cB^{*}, s^{*}) \in \arg\max\limits_{\substack{\cA, \cB \subseteq \Omega, s \in \Omega \\|\cA| = \ell, |\cB| = k}} \Big(\alpha f(\cA \cup \cB \cup \{s\}) - \alpha f(\cA \cup \cB) - \alpha f(\cA\cup\{s\}) + \alpha f(\cA)\Big)$, which implies $\metric(\alpha f) = d^{\ell, k}(\alpha f) = \alpha \metric[f]$, which proves positive homogeneity.

Now suppose $\metric[f] > \epsilon$, for some $\epsilon \geq 0$. Let $g:2^{\Omega}\to \mathbb{R}$ be any nonnegative, increasing set function such that $||g - f||_{\infty} < \epsilon/4$. Consider $(\cA^{*}, \cB^{*}, s^{*}) \in \arg\max\limits_{\substack{\cA, \cB \subseteq \Omega, s \in \Omega \\|\cA| = \ell, |\cB| = k}} \Big(f(\cA \cup \cB \cup \{s\}) - f(\cA \cup \cB) - f(\cA\cup\{s\}) + f(\cA)\Big)$, 
we have
\begin{align*}
\metric(g) &\geq g(\cA^{*} \cup \cB^{*} \cup\{s^{*}\}) - g(\cA^{*} \cup \cB^{*}) - g(\cA^{*} \cup\{s^{*}\}) + g(\cA^{*})\\
&\geq f(\cA^{*} \cup \cB^{*} \cup\{s^{*}\}) - f(\cA^{*} \cup \cB^{*}) - f(\cA^{*} \cup\{s^{*}\}) + f(\cA^{*}) - \epsilon 
= \metric[f] - \epsilon 
> 0,
\end{align*}
which implies $g$ is not submodular.
$\myQED$

\begin{makeshiftResult}{induceASN}
Define $\glob_{+}:\mathcal{F} \to \mathbb{R}$ by $\glob_{+}[f] \coloneqq \max\{0, \glob[f]\}$. Then $\glob_{+}$ is an asymmetric seminorm on $\mathcal{F}$. 
\end{makeshiftResult}

\myProof
Let $f_{1}, f_{2} \in \mathcal{F}$. Observe that $\glob_{+}[f_{1}] + \glob_{+}[f_{2}] \geq \glob[f_{1}] + \glob[f_{2}] \geq \glob[f_{1} + f_{2}]$, by \thref{ApproxSubmodSubadditive}. If $\glob[f_{1} + f_{2}] \geq 0$, then $\glob_{+}[f_{1}] + \glob_{+}[f_{2}] \geq \glob[f_{1} + f_{2}] = \glob_{+}[f_{1} + f_{2}]$. Otherwise, $\glob_{+}[f_{1} + f_{2}] = 0$. We have $\glob_{+}[f_{1}] \geq 0$, $\glob_{+}[f_{2}] \geq 0$, which imply $\glob_{+}[f_{1}] + \glob_{+}[f_{2}] \geq 0 = \glob_{+}[f_{1} + f_{2}]$, thus proving subadditivity.

Suppose $f \in \mathcal{F}$ is such that $\glob[f] \geq 0$. By \thref{ApproxSubmodSubadditive}, for any $\alpha \geq 0, \alpha\glob_{+}[f] = \alpha\glob[f] = \glob[\alpha f] = \glob_{+}[\alpha f]$. If $\glob[f] < 0$, then $\glob_{+}[f] = \glob[f] = 0.$ By similar arguments to those of \thref{ApproxSubmodSubadditive}, $\alpha\glob_{+}[f] = 0 = \glob_{+}[\alpha f]$; hence $\glob_{+}$ satisfies positive homogeneity. Clearly, $\glob_{+}$ is nonnegative, which concludes the proof.
$\myQED$

\begin{lemma}{\emph{(Narayanan 1997 \citep{Narayanan1997})}} \thlabel{modularSetProperty}
For any modular set function $g':2^{\Omega} \to \mathbb{R}$ and $\cA, \cB, \cS, \cT \subseteq \Omega$ with $\cA \cup \cB = \cS \cup \cT$ and $\cA \cap \cB = \cS \cap \cT$, we have $g'(\cA) + g'(\cB) = g'(\cS) + g'(\cT)$.
\end{lemma}

\begin{repeattheorem}{submodPreserve}
Given $f:2^{\Omega} \to \mathbb{R}$ and the corresponding functions $f_{1}, f_{2}, f_{\cA},$ and $f_{q}$,
we have:
\begin{enumerate*}[label=(\roman*)]
\item $\glob[f] = \glob[f_{1}]$. 
\item $2\glob[f] \geq \glob[f_{2}]$. 
\item $\glob[f] \geq \glob[f_{\cA}]$ 
\item $\glob[f] \geq \glob[f_{q}]$. 
\item $\glob[f] \geq \glob[f\circledast g]$. 
\end{enumerate*}
\end{repeattheorem}

\myProof
\ref{submodPreserve1}: For any $\cA, \cB \subseteq \Omega$, by de Morgan's laws, $(\cA \cup \cB)^{c} = \cA^{c} \cap \cB^{c}$ and $(\cA \cap \cB)^{c} = \cA^{c} \cup \cB^{c}$, from which the result immediately follows. 

\ref{submodPreserve2}: Define $\tilde{f}: 2^{\Omega} \to \mathbb{R}$ where $\tilde{f}(\cS) = f(\cS^{c}) - f(\Omega)$. Because $f(\Omega)$ is a constant, by \ref{submodPreserve1}, $\glob[\tilde{f}] = \glob[f]$. Thus, by \thref{ApproxSubmodSubadditive}, $\glob[f_{2}] \leq \glob[f] + \glob[\tilde{f}] = 2\glob[f]$.

\ref{submodPreserve3}: Let $\cS, \cT \subseteq \cA^{c}$. Then $\cS \cup \cA$ and $\cT \cup \cA$ are subsets of $\Omega$. Hence, $f_{\cA}(\cS) + f_{\cA}(\cT) + \glob[f] = f(\cA \cup \cS) + f(\cA \cup \cT) + \glob[f] \geq f((\cA \cup \cS) \cup (\cA \cup \cT)) + f((\cA \cup \cS) \cap (\cA \cup \cT)) = f(\cA \cup (\cS \cup \cT)) + f(\cA \cup (\cS \cap \cT)) = f_{\cA}(\cS \cup \cT) + f_{\cA}(\cS \cap \cT)$.

\ref{submodPreserve4}: Consider $\cA, \cB \subseteq \Omega(q)$. We have $\glob[f] + f_{q}(\cA) + f_{q}(\cB) = \glob[f] + f(\bigcup\limits_{i \in \cA} \cS(i)) + f(\bigcup\limits_{i \in \cB} \cS(i)) \leq f(\bigcup\limits_{i \in \cA \cup \cB} \cS(i)) + f(\bigcup\limits_{i \in \cA \cap \cB} \cS(i)) = f_{q}(\cA \cup \cB) + f_{q}(\cA \cap \cB).$

\ref{submodPreserve5} This proof is similar to that of \cite{Narayanan1997}. Let $\cS, \cT, \cZ_{\cS}, \cZ_{\cT} \subseteq \Omega,$ where $f \circledast g(\cS) = f(\cZ_{\cS}) + g(\cC \backslash \cZ_{\cS}),$ and $f \circledast g(\cT) = f(\cZ_{\cT}) + g(\cT \backslash \cZ_{\cT})$. It is not hard to show that $(\cS \backslash \cZ_{\cS}) \cup (\cT \backslash \cZ_{\cT}) = ((\cS \cup \cT) \backslash (\cZ_{\cS} \cup \cZ_{\cT})) \cup ((\cS \cap \cT) \backslash ((\cZ_{\cS} \cap \cZ_{\cT}))$ and $(\cS \backslash \cZ_{\cS})\cap (\cT\backslash\cZ_{\cT}) = ((\cS \cup \cT)\backslash (\cZ_{\cS} \cup \cZ_{\cT})) \cap ((\cS \cap \cT) \backslash (\cZ_{\cS} \cap \cZ_{\cT})$.

By \thref{modularSetProperty}, 
$g(\cS \backslash \cZ_{\cS}) + g(\cT \backslash \cZ_{\cT}) = g((\cS \cup \cT)\backslash(\cZ_{\cS}\cup \cZ_{\cT})) + g((\cS \cap \cT)\backslash(\cZ_{\cS}\cap \cZ_{\cT})).$

By the definition of $\cZ_{\cS}$ and $\cZ_{\cT}$, \begin{align*}
&f\circledast g(\cS) + f\circledast g(\cT)\\
= \ &f(\cZ_{\cS}) + f(\cZ_{\cT}) + g((\cS \cup \cT)\backslash(\cZ_{\cS}\cup \cZ_{\cT})) + g((\cS \cap \cT)\backslash(\cZ_{\cS}\cap \cZ_{\cT}))\\
\geq &-\glob[f] + f(\cZ_{\cS} \cup \cZ_{\cT}) + f(\cZ_{\cS} \cap \cZ_{\cT}) + g((\cS \cup \cT)\backslash(\cZ_{\cS}\cup \cZ_{\cT})) + g((\cS \cap \cT)\backslash(\cZ_{\cS}\cap \cZ_{\cT}))\\
\geq &-\glob[f] + \min\limits_{\cZ \subseteq \cS \cup \cT} f(\cZ) + g((\cS \cup \cT) \backslash \cZ) + \min\limits_{\cZ \subseteq \cS \cap \cT} f(\cZ) + g((\cS \cap \cT) \backslash \cZ)\\
= &-\glob[f] + f \circledast g(\cS \cup \cT) + f \circledast g(\cS \cap \cT).
\end{align*}
Hence, $\glob[f] \geq \glob[f \circledast g]$.
$\myQED$

\begin{lemma}\emph{{(Bach 2013 \citep{Bach2013}})}\thlabel{LovaszHomogeneous}
Given a set function $f$, its Lov\'{a}sz extension is positively homogeneous of degree 1.
\end{lemma}

\begin{lemma}\thlabel{subtractOffSet}
Suppose $x = w + \alpha x(\cA) \in [0,1]^{|\Omega|}$, where $w, \alpha x(\cA) \in [0,1]^{|\Omega|}, \cA \subseteq \Omega$, $\alpha \in \mathbb{R}_{+}$, and $x_{s} \geq x_{s'}$ for any $s \in \cA, s' \in \cA^{c}$. If there exists a permutation $\pi$ of $\Omega$ such that $x_{\pi_{1}} \geq \cdots \geq x_{\pi_{|\Omega|}}$ and $w_{\pi_{1}} \geq \cdots \geq w_{\pi_{|\Omega|}}$, then $F^{L}(x) = F^{L}(w) + \alpha f(\cA)$.
\end{lemma}

\myProof
Let $\pi$ be the ranking permutation in the hypothesis, and note that it is also a ranking permutation of $x(\cA)$; i.e., $x(\cA)_{\pi_{1}} \geq \cdots \geq x(\cA)_{\pi_{|\Omega|}}$. We have $F^{L}(x) = \sum\limits_{k = 1}^{|\Omega|} x_{\pi_{k}}(f(x(\cS^{\pi}_{k})) - f(x(\cS^{\pi}_{k-1})))$, and $F^{L}(w) = \sum\limits_{k = 1}^{|\Omega|} w_{\pi_{k}}(f(x(\cS^{\pi}_{k})) - f(x(\cS^{\pi}_{k-1})))$, implying that $F^{L}(x) - F^{L}(w) = \sum\limits_{k = 1}^{|\Omega|} \alpha x(\cA)_{\pi_{k}}(f(x(\cS^{\pi}_{k})) - f(x(\cS^{\pi}_{k-1})) = F^{L}(\alpha x(\cA)) = \alpha F^{L}(x(\cA)) = \alpha f(\cA),$ where we have used the positive homogeneity of $F^{L}$ (\thref{LovaszHomogeneous}).
$\myQED$

\begin{lemma}\thlabel{inductionGreaterThan}
Let $f:2^{\Omega} \to \mathbb{R}$ be increasing with $f(\emptyset) = 0$, and let $\tilde{\gamma} \in \Gamma(f)$. For any $\cS \subseteq \Omega$, $f(\cS) \geq -|\cS|\dknap[f] + \sum\limits_{s \in \cS}\tilde{\gamma}_{s}$.
\end{lemma}

\myProof
Consider a permutation $(\rho_{1},\dots,\rho_{|\Omega|})$ such that $\tilde{\gamma}_{\rho_{1}}\geq \cdots \geq \gamma_{\rho_{|\Omega|}}$ and set $\theta^{*} = -|\Omega|\dknap[f]$. We prove by induction on $|\cS|$ that $f(\cS) \geq -|\cS|\dknap[f] + \sum\limits_{s \in \cS} \tilde{\gamma}_{s},$ for all $\cS \subseteq \Omega$. The base case is confirmed as $f(\emptyset) = 0 = \sum\limits_{s \in \emptyset}\tilde{\gamma}_{s}$. Assume for all $\tilde{\cS} \subseteq \Omega$ with $|\tilde{\cS}| \leq \alpha$, $f(\tilde{\cS}) \geq -|\tilde{\cS}|\dknap[f] + \sum\limits_{s \in \tilde{\cS}} \tilde{\gamma}_{s}$, and let $|\cS| = \alpha + 1$. Set $k = \max\{i \setbar \rho_{i} \in \cS\}$. Then $\cS \cup \cS^{\rho}_{k - 1} = \cS^{\rho}_{k}$ and $\cS\cap \cS^{\rho}_{k-1} = \cS\backslash \{\rho_{k}\}$. Observe that by the definition of $\dknap[f]$,
\begin{align*}
f(\cS) &\geq f(\cS \cup \cS^{\rho}_{k-1}) + f(\cS \cap \cS^{\rho}_{k-1}) - f(\cS^{\rho}_{k-1}) - \dknap[f]\\
&= f(\cS^{\rho}_{k}) - f(\cS^{\rho}_{k-1}) - \dknap[f] + f(\cS \backslash \{\rho_{k}\})\\
&= \tilde{\gamma}_{\rho_{k}} + f(\cS\backslash\{\rho_{k}\}) - \dknap[f]\\
&\geq \sum\limits_{s \in \cS}\tilde{\gamma}_{s} - |\cS|\dknap[f],
\end{align*}
where the last line uses the induction hypothesis. $\myQED$

\begin{repeattheorem}{LovaszApprox}
For any increasing set function $f:2^{\Omega} \to \mathbb{R}$ such that $f(\emptyset) = 0$, 
\begin{align*}
F^{L}(x) &\leq \max\limits_{\gamma \in \Gamma(f)} \sum\limits_{s \in \Omega} \gamma_{s}x_{s}
\leq F^{C}(x) + |\Omega|\dknap[f]\\
&\leq F^{L}(x) + |\Omega|\dknap[f] 
\leq \max\limits_{\gamma \in \Gamma(f)} \sum\limits_{s \in \Omega} \gamma_{s}x_{s} + |\Omega|\dknap[f].
\end{align*}
Hence, $F^{L}(x) \geq F^{C}(x) \geq F^{L}(x) - |\Omega|\dknap[f]$, and $||F^{L} - F^{C}||_{\infty} \leq |\Omega|\dknap[f]$. Moreover, $F^{L}$ is $|\Omega|\dknap[f]$-approximately convex.

In addition, if for some $\epsilon > 0$, $F^{L}$ is $\epsilon$-approximately convex,  then $\dknap[f] \leq \epsilon$.
\end{repeattheorem}

\myProof
The following proof uses a version of well-known linear programming duality arguments (e.g., \cite{Lovasz1983}). We first suppose $f$ is not submodular (hence $\dknap[f] > 0$).

Given $x \in [0,1]^{|\Omega|}$ there exists a permutation $(\pi_{1},\dots,\pi_{|\Omega|})$ such that $x_{\pi_{1}} \geq \dots\geq x_{\pi_{|\Omega|}}$. Let $x_{\pi_{0}} = 1$. Consider the dual of \eqref{closureLPDual}:
\begin{align}
\max\limits_{\gamma,\theta}\left\{\theta + \sum\limits_{s \in \Omega} x_{s}\gamma_{s} \bigg\vert \theta + \sum\limits_{s \in \cS} x_{s} \leq f(\cS), \forAll \cS \subseteq \Omega \right\}. \label{closureLPPrimal}
\end{align}
Define $y^{*} \in \mathbb{R}^{|2^{\Omega}|}$ by 
$y^{*}_{\cS}$ equals $x_{\pi_{i}} - x_{\pi_{i + 1}}, \ \text{if } \cS = \cS^{\pi}_{i}, i \in \{0,\dots, |\Omega| - 1\}$, $x_{\pi_{|\Omega|}}, \ \text{if } \cS = \Omega$, and 0 otherwise.
\cutOut{\begin{align*}
y^{*}_{\cS} &= \begin{cases}
&x_{\pi_{i}} - x_{\pi_{i + 1}}, \ \text{if } \cS = \cS^{\pi}_{i}, i \in \{0,\dots, |\Omega| - 1\},\\
&x_{\pi_{|\Omega|}}, \ \text{if } \cS = \Omega,\\
&0, \ \text{otherwise}. 
\end{cases}
\end{align*}
}

We first show $y^{*}$ is feasible for \eqref{closureLPDual}. 
Observe that $\sum\limits_{s \in \Omega}y^{*}(\cS) = \sum\limits_{i = 0}^{|\Omega| - 1}(x_{\pi_{i}} - x_{\pi_{i+1}}) + x_{\pi_{|\Omega|}} = x_{\pi_{0}} = 1$.
In addition, for any $s \in \Omega, s = \pi_{j}$ for some $j \in \Omega;$ hence, $\sum\limits_{\cS \ni s}y^{*}(\cS) = \sum\limits_{i = j}^{|\Omega|} y^{*}(\cS^{\pi}_{i}) = x_{\pi_{j}} = x_{s}$. 
Moreover, it is easy to observe that $y^{*}$ is nonnegative. Hence, $y^{*}$ is feasible for \eqref{closureLPDual}, and $F^{C}(x) \leq \sum\limits_{\cS \subseteq \Omega}f(\cS)y^{*}(\cS) = x_{\pi_{|\Omega|}}f(\Omega) + \sum\limits_{i = 0}^{|\Omega| - 1}(x_{\pi_{i}} - x_{\pi_{i+1}})f(\cS^{\pi}_{i}) = F^{L}(x)$.

Consider $\gamma^{*}\in \Gamma(f)$ such that $\gamma^{*}_{\pi_{1}} \geq \cdots\geq \gamma^{*}_{\pi_{|\Omega|}}$. Then $\sum\limits_{s \in \Omega}\gamma^{*}_{s}x_{s} = \sum\limits_{ i = 1}^{|\Omega|}(f(\cS^{\pi}_{i}) - f(\cS^{\pi}_{i-1}))x_{\pi_{i}} = F^{L}(x)$. By \thref{inductionGreaterThan}, $(\gamma^{*}, -|\Omega|\dknap[f])$ is feasible for \eqref{closureLPPrimal}.

Therefore, 
$F^{L}(x) \leq \max\limits_{\gamma \in \Gamma(f)} \ \sum\limits_{s \in \Omega} \gamma_{s}x_{s} 
\leq F^{C}(x) + |\Omega|\dknap[f] 
\leq F^{L}(x) + |\Omega|\dknap[f] 
= |\Omega|\dknap[f] + \max\limits_{\gamma \in \Gamma(f)} \ \sum\limits_{s \in \Omega} \gamma_{s}x_{s}.$
This also implies that $F^{L}(x) \geq F^{C}(x) \geq F^{L}(x) - |\Omega|\dknap[f]$ and $||F^{L} - F^{C}||_{\infty} \leq |\Omega|\dknap[f]$.

To show that $F^{L}$ is approximately convex, consider $x, y \in [0,1]^{|\Omega|}, \lambda \in [0,1]$, then we have
\begin{align*}
F^{L}(\lambda x + (1 - \lambda)y) &\leq F^{C}(\lambda x + (1 - \lambda)y) + |\Omega|\dknap[f]\\
&\leq \lambda F^{C}(x) + (1 - \lambda)F^{C}(y) + |\Omega|\dknap[f]\\
&\leq \lambda F^{L}(x) + (1 - \lambda)F^{L}(y) + |\Omega|\dknap[f].
\end{align*}

For the last statement, suppose $F^{L}$ is $\epsilon$-approximately convex. For any $\cA, \cB \in \Omega$, denote the symmetric difference as $\cA \dotminus \cB = (\cA \backslash \cB) \cup (\cB \backslash \cA)$. Consider any $S, T \subseteq \Omega$, then we have 
\begin{align*}
\frac{1}{2}(f(\cS) + f(\cT)) + \epsilon &= \frac{1}{2}F^{L}(x(\cS)) + \frac{1}{2}F^{L}(x(\cT)) + \epsilon \\
&\geq F^{L}(\frac{1}{2}(x(\cS) + x(\cT))) \\
&= F^{L}(\frac{1}{2}(x(\cS \cap \cT) + x(\cS \cap \cT) + x(\cS \dotminus \cT)))\\
&= F^{L}(\frac{1}{2}(x(\cS \cap \cT) + x(\cS \cup \cT)))\\
&= \frac{1}{2}(f(\cS \cap \cT) + f(\cS \cup \cT))
\end{align*}
where the last line uses the fact that the Lov\'{a}sz extension is positively homogeneous \thref{LovaszHomogeneous} and \thref{subtractOffSet}. 

Suppose $f$ is submodular. Then, $\dknap[f] = 0$ and $F^{C} = F^{L}$ (e.g., see \cite{Lovasz1983,Bach2013}). A similar linear programming duality argument (with $\dknap[f]$ replaced with 0) proves that $F^{L}(x) = \max\limits_{\gamma \in \Gamma(f)} \sum\limits_{s \in \Omega} \gamma_{s}x_{s}.$
$\myQED$


\begin{repeattheorem}{SubmodToConvexLovasz}
Given set functions $f,g:2^{\Omega} \to \mathbb{R}$, where $f(\emptyset) = g(\emptyset) = 0$, and their respective Lov\'{a}sz extensions $F^{L},G^{L}: [0,1]^{|\Omega|} \to \mathbb{R}$, $||F^{L} - G^{L}||_{\infty} = ||f - g||_{\infty}$.
\end{repeattheorem}
\myProof
Let $x \in [0,1]^{|\Omega|}$ and $\pi$ a permutation of $(1,\dots,|\Omega|)$ such that $x_{\pi_{1}} \geq x_{\pi_{2}} \geq \dots \geq x_{\pi_{|\Omega|}},$ and set $x_{\pi_{|\Omega|+1}} = 0$. Then $F^{L}(x) = \sum\limits_{i = 1}^{|\Omega|} (f(\cS^{\pi}_{i}) - f(\cS^{\pi}_{i - 1}))x_{\pi_{i}} = -x_{\pi_{|\Omega|}}f(\Omega) + \sum\limits_{i = 1}^{|\Omega| - 1} (x_{\pi_{i}} - x_{\pi_{i+1}})f(\cS^{\pi}_{i}),$ and similarly for $G^{L}$. Hence, 
\begin{align*}
|F^{L}(x) - G^{L}(x)| &= \left|x_{\pi_{|\Omega|}}(f(\Omega) - g(\Omega)) + \sum\limits_{i = 1}^{|\Omega|-1} (x_{\pi_{i}} - x_{\pi_{i+1}})(f(\cS^{\pi}_{i}) - g(\cS^{\pi}_{i}))\right|\\
&\leq |x_{\pi_{|\Omega|}}(f(\Omega) - g(\Omega))| + \sum\limits_{i = 1}^{|\Omega|-1} \left|(x_{\pi_{i}} - x_{\pi_{i+1}})(f(\cS^{\pi}_{i}) - g(\cS^{\pi}_{i}))\right|\\
&= |x_{\pi_{|\Omega|}}|\cdot|(f(\Omega) - g(\Omega))| + \sum\limits_{i = 1}^{|\Omega|-1} |(x_{\pi_{i}} - x_{\pi_{i+1}})|\cdot|(f(\cS^{\pi}_{i}) - g(\cS^{\pi}_{i}))|\\
&\leq |x_{\pi_{|\Omega|}}|\cdot||f - g||_{\infty} + \sum\limits_{i = 1}^{|\Omega|-1} |(x_{\pi_{i}} - x_{\pi_{i+1}})|\cdot||f - g||_{\infty}\\
&= ||f - g||_{\infty}\left(x_{\pi_{|\Omega|}} + \sum\limits_{i = 1}^{|\Omega|-1} (x_{\pi_{i}} - x_{\pi_{i+1}})\right)\\
&\leq ||f - g||_{\infty}.
\end{align*}
Hence, $||F^{L} - G^{L}||_{\infty} \leq ||f - g||_{\infty}$. Moreover, for some $\cS \subseteq \Omega, ||f - g||_{\infty} = |f(\cS) - g(\cS)| = |F^{L}(x(\cS)) - G^{L}(x(\cS))|$, which implies $||F^{L} - G^{L}||_{\infty} = ||f - g||_{\infty}.$ $\myQED$

\begin{repeattheorem}{generalInequality}
For any $\gamma \in \Gamma(g_{\tau}), \cS \subseteq \Omega,$ we have $-|\Omega|\dknap[f_{\tau}] + \sum\limits_{s \in \Omega} \gamma_{s}(x(\cS))_{s} \leq f_{\tau}(\cS) - \phi(\tau)$.
\end{repeattheorem}

\myProof
Recall that $g_{\tau}(\emptyset) = 0$ and $\dknap[f_{\tau}] = \dknap[g_{\tau}]$. By \thref{LovaszApprox}, $-|\Omega|\dknap[f_{\tau}] + \sum\limits_{s \in \Omega} \gamma_{s}(x(\cS))_{s} \leq G^{L}_{\tau}(x(\cS)) = g_{\tau}(\cS) = f_{\tau}(\cS) - \phi(\tau).$
$\myQED$

\begin{repeattheorem}{validEpigraphBinary}
For any $\gamma \in \Gamma(g_{\sigma})$, the following inequality is valid for $\conv(H_{\mathbb{B}})$:
\begin{align}\label{validSigmaInequality}
-|\Omega|\dknap[f_{\sigma}] + \sum\limits_{s \in \Omega}\gamma_{s}x_{s} \leq z - \phi(\sigma).
\end{align} 
\end{repeattheorem}

\myProof
This proof follows arguments similar to that of \cite{Atamturk2008}. Consider $(x,z) \in H_{\mathbb{B}}$, which implies $x = x(\cS)$, for some $\cS \subseteq \Omega$. From \thref{generalInequality}, $-|\Omega|\dknap[f_{\sigma}] + \sum\limits_{s \in \Omega}\gamma_{s}(x(\cS))_{s} \leq f_{\sigma}(\cS) - \phi(\sigma) = g_{\sigma}(\cS)$. Because $(x,z) \in H_{\mathbb{B}}, F_{\sigma}(x) = f_{\sigma}(\cS) = g_{\sigma}(\cS) + \phi(\sigma) \leq z$, which implies $-|\Omega|\dknap[f_{\sigma}] + \sum\limits_{s \in \Omega}\gamma_{s}x_{s} \leq z - \phi(\sigma).$
$\myQED$

\begin{repeattheorem}{polarIneq}
Let $f:2^{\Omega} \to \mathbb{R}$ be increasing with $f(\emptyset) = 0$. Suppose $\gamma \in \mathbb{R}^{|\Omega|}$ and
\begin{align}
\sum\limits_{s \in \Omega}\gamma_{s}x_{s} \leq z + |\Omega|\dknap[f]+\gamma_{0} \label{ineqGeneralOut}
\end{align}
defines a nontrivial facet of $H_{f}$. Let $\bar{f}:2^{\Omega} \to \mathbb{R}$ be defined by $\bar{f}(\emptyset) = 0, \bar{f}(\cS) = f(\cS) + |\Omega|\dknap[f] + \gamma_{0}$, for all nonempty $\cS \subseteq \Omega$, and suppose $\gamma \in \Gamma(\bar{f})$. Then, $\gamma_{0} \leq 0$.
\end{repeattheorem}

\myProof
This proof follows steps similar to that of \cite{AtamturkShort2020}, with additional arguments to account for approximate submodularity. Suppose $\gamma_{0} > 0$. Because \eqref{ineqGeneralOut} is a valid inequality for $H_{f}$, for any non-empty $\cS \subseteq \Omega,\sum\limits_{s \in \cS} \gamma_{s} = \sum\limits_{s \in \Omega}\gamma_{s}x(\cS)_{s} \leq f(\cS) + \gamma_{0} + |\Omega|\dknap[f] = \bar{f}(\cS),$ and $\sum\limits_{s \in \emptyset} \gamma_{s} = 0 = \bar{f}(\emptyset)$. Thus, $\gamma \in P_{\bar{f}}$, and by \thref{atamturkOuter}, \eqref{validIneqOne} is valid for $H_{\bar{f}}$:
\begin{align}
\sum\limits_{s \in \Omega}\gamma_{s}x_{s} \leq z. \label{validIneqOne}
\end{align}

We show that \eqref{validIneqOne} is facet-defining for $H_{\bar{f}}$. Observe that the solutions $\{(x(\{s\}), f(\{s\}))\}_{s \in \Omega}$ $\cup \{(x(\emptyset), 1),$ $(x(\emptyset),0)\}$ are $|\Omega| + 2$ affinely independent solutions, so the dimension of $H_{f}$ is $|\Omega|+1$. By the hypothesis, $\sum\limits_{s \in \Omega}\gamma_{s}x_{s} \leq z + |\Omega|\dknap[f] + \gamma_{0}$ is facet-defining for $H_{f}$. Thus, there exist $|\Omega| + 1$ affinely independent solutions $\{(x^{k},z^{k})\}_{k = 1}^{|\Omega| + 1}$ such that $(x^{k},z^{k}) \in H_{f}$ and $\sum\limits_{s \in \Omega}\gamma_{s}x^{k}_{s} = z^{k} + |\Omega|\dknap[f] + \gamma_{0}$. By Carath\'{e}odory's theorem, each of the affinely independent solutions $(x^{k},z^{k})$ can be represented by a convex combination of $|\Omega| + 2$ (integral) extreme points of $H_{f}$: $(x^{k}, z^{k}) = \sum\limits_{\ell = 1}^{|\Omega| + 2} \lambda^{k,\ell}(x^{k,\ell},z^{k,\ell})$, where $\sum\limits_{\ell = 1}^{|\Omega| + 2}\lambda^{k,\ell} = 1$, $\lambda^{k} \in \mathbb{R}^{|\Omega| + 2}_{+},$ and $(x^{k,\ell}, z^{k,\ell})$ is an integral extreme point of $H_{f}$, for all $\ell \in \{1,\dots,|\Omega| + 2\}$. 
Consider $(x^{k,\ell},z^{k,\ell})$ for some $k \in \{1,\dots,|\Omega|+1\}$ and $\ell \in \{1,\dots,|\Omega| + 2\}$. Suppose that $x^{k,\ell} = 0$; because $(x^{k,\ell},z^{k,\ell}) \in H_{f}$, $0 \leq z^{k,\ell}$. Thus, $\bar{f}(x^{k,\ell}) = 0 \leq z^{k,\ell} < z^{k,\ell} + |\Omega|\dknap[f] + \gamma_{0}$. If instead $x^{k,\ell} \neq 0,$ then $x^{k,\ell} = x(\cS)$ for some nonempty $\cS \subseteq \Omega$. Notice that by $(x^{k,\ell},z^{k,\ell}) \in H_{f}$, $f(x^{k,\ell}) \leq z^{k,\ell}$, so $\bar{f}(x^{k,\ell}) \leq z^{k,\ell} + |\Omega|\dknap[f] + \gamma_{0}$. Hence, $(x^{k,\ell},z^{k,\ell} + |\Omega|\dknap[f] + \gamma_{0}) \in H_{\bar{f}}$; moreover, $(x^{k},z^{k} + |\Omega|\dknap[f] + \gamma_{0}) = \sum\limits_{\ell = 1}^{|\Omega| + 2} \lambda^{k,\ell}(x^{k,\ell},z^{k,\ell} + |\Omega|\dknap[f] + \gamma_{0}) \in H_{\bar{f}}$.

Suppose that the points $(x^{k}, z^{k} + |\Omega|\dknap[f] + \gamma_{0})_{k = 1}^{|\Omega| + 1}$ are not affinely independent. Then there exists $\sigma \in \mathbb{R}^{|\Omega| + 1} \backslash \{0\}$ such that $\sum\limits_{k = 1}^{|\Omega| + 1}\sigma_{k} = 0$ and $\sum\limits_{k = 1}^{|\Omega| + 1} \sigma_{k}(x^{k}, z^{k}  + |\Omega|\dknap[f] + \gamma_{0}) = (0,0)$. Let $j \in \{1,\dots,|\Omega| + 1\}$ be such that $\sigma_{j} \neq 0$; without loss of generality, let $\sigma_{j} = 1$. Thus, $\sum\limits_{k \neq j}\sigma_{k}(x^{k}, z^{k}  + |\Omega|\dknap[f] + \gamma_{0}) = -(x^{j},z^{j}  + |\Omega|\dknap[f] + \gamma_{0}).$ Because $\sigma_{j} = 1$, $\sum\limits_{k \neq j} \sigma_{k} = -1$; thus, 
\begin{align*}
-(x^{j}, z^{j} + |\Omega|\dknap[f] + \gamma_{0}) &= \sum\limits_{k \neq j}\sigma_{k}(x^{k}, z^{k}  + |\Omega|\dknap[f] + \gamma_{0}) \\
&= -(0, |\Omega|\dknap[f] + \gamma_{0}) + \sum\limits_{k \neq j}\sigma_{k}(x^{k},z^{k})\\
\iff -(x^{j}, z^{j}) &= \sum\limits_{k \neq j} \sigma_{k}(x^{k},z^{k}),
\end{align*}
which contradicts the affine independence of $(x^{k},z^{k})_{k = 1}^{|\Omega| + 1}$, so $(x^{k}, z^{k} + |\Omega|\dknap[f] + \gamma_{0})_{k = 1}^{|\Omega| + 1}$ are affinely independent.
Also, $\sum\limits_{s \in \Omega}{\gamma}_{s}x^{k}_{s} = z^{k} + |\Omega|\dknap[f] +\gamma_{0}$, for each $k \in \{1,\dots,|\Omega| + 1\}$, which implies \eqref{validIneqOne} is facet-defining for $H_{\bar{f}}$. 

By \thref{atamturkOuter}, $\gamma$ is an extreme point of $P_{\bar{f}}$, and by the hypothesis, there exists a permutation $(\rho_{1},\dots,\rho_{|\Omega|})$ such that $\gamma_{\rho_{s}}=\bar{f}(\cS^{\rho}_{s}) - \bar{f}(\cS^{\rho}_{s - 1})$, for all $s \in \Omega$. Define $\hat{\gamma}$ by $\hat{\gamma}_{\rho_{1}} = \gamma_{\rho_{1}} - \gamma_{0}, \hat{\gamma}_{\rho_{s}} = \gamma_{\rho_{s}}$, otherwise. Thus, $\hat{\gamma}_{\rho_{s}} = f(\cS^{\rho}_{s}) - f(\cS^{\rho}_{s-1})$ and $\hat{\gamma} \in \Gamma(f)$. Because $\gamma_{0} > 0$, $\bar{f}$ is increasing, so that by \thref{inductionGreaterThan}, $\sum\limits_{s \in \cS}\hat{\gamma}_{s} \leq f(\cS) + |\cS|\dknap[f]$. Hence, if $\gamma'_{s} = \hat{\gamma}_{s} - \dknap[f]$, then $\gamma' \in P_{f}$. By \thref{atamturkOuter}, we have
$\sum\limits_{s \in \Omega} \gamma'_{s}x_{s} \leq z$
is valid for $H_{f}$, which implies $(\gamma_{\rho_{1}} - \gamma_{0})x_{\rho_{1}} + \sum\limits_{s \in \Omega \backslash \{\rho_{1}\}} \hat{\gamma}_{s}x_{s} \leq z + |\Omega|\dknap[f]$ is also valid for $H_{f}$. Because $\gamma_{0} > 0$, we also have the valid inequality $\gamma_{0}x_{\rho_{1}} \leq \gamma_{0}$.

Combining these last two inequalities implies 
$\sum\limits_{s \in \Omega} \gamma_{s}x_{s} \leq z + |\Omega|\dknap[f] + \gamma_{0},$
thus the facet-defining inequality \eqref{ineqGeneralOut} is dominated, a contradiction.
$\myQED$

\begin{proposition}\thlabel{XFullDimensional}
If $f(\{s\}) \leq b,$ for all $s \in \Omega$, then $X$ is full-dimensional.
\end{proposition}

\myProof
By the hypothesis, the zero vector and $x(\{s\})$ are feasible for each $s \in \Omega$. Hence there are $|\Omega| + 1$ affinely independent points in $X$, implying the dimension of $X$ is $|\Omega|$.
$\myQED$

\begin{proposition}\thlabel{thm:easyInequalities}\emph{{(Hammer et al. 1975 \citep{Hammer1975}, Atamt\"{u}rk and Narayanan 2009 \citep{Atamturk2009})}}
\begin{enumerate}
\item The inequality $x(\{s\}) \geq 0$ is facet-defining for $\conv(X)$, for all $s \in \Omega$.
\item The inequality $x(\{s\}) \leq 1$ is facet-defining for $\conv(X)$ if and only if $f(\{s,t\}) \leq b$ for all $t \in \Omega \backslash s$.
\end{enumerate}
\end{proposition}

\begin{repeattheorem}{Atamturk09ExtendedV2}
If $\cS \subseteq \Omega$ is a cover for $X$, the extended cover inequality $\sum\limits_{s \in E_{\pi}(\cS)} x_{s} \leq |\cS| - 1$ is valid for $X$ if $f(\cS) > (|\cS| + |U_{\pi}(\cS)|)\dknap[f] + b$. In addition, the inequality defines a facet of $\{x \in X \setbar x_{s} = 0, \forAll s \not\in E_{\pi}(\cS)\}$ if $\cS$ is also a minimal cover and for each $s \in U_{\pi}(\cS)$, there exist $t_{s}, u_{s} \in \cS$ such that $t_{s} \neq u_{s}$, and $f(\cS \cup \{s\} \backslash \{t_{s},u_{s}\}) \leq b$.
\end{repeattheorem}

\myProof
This proof uses steps similar to those in \cite{Atamturk2009}, who establish the submodular case. We show that if $x \in [0,1]^{|\Omega|}$ with $\sum\limits_{s \in E_{\pi}(\cS)} x_{s} > |S| - 1$, then $x_{s} \not\in X$. Because $X$ is the convex hull of characteristic vectors, it suffices to consider such characteristic vectors. That is $x(\tilde{S})$, where $\tilde{S} \subseteq \Omega$ and there exists $\cT \subseteq \tilde{S}$ such that $\cT \subseteq E_{\pi}(\cS)$ with $|\cT| \geq |\cS|$. In this case, $\sum\limits_{s \in E_{\pi}(\cS)} x(\tilde{S})_{s} \geq \sum\limits_{s \in \cT} x(\tilde{S})_{s} \geq |\cS|$. Let $\cK = \cS \backslash \cT$, and $\cL = U_{\pi}(\cS) \cap \cT = \{\ell_{1},\dots,\ell_{|\cL|}\}$, with indexing consistent with $\pi$. 

Observe that $\cS \backslash \cK = \cS \cap \cT$ and $(\cS \cup \cL) \backslash \cK = (\cS \cup (U_{\pi}(\cS) \cap \cT)) \backslash (\cS \backslash \cT) = (\cS \cap \cT) \cup (U_{\pi}(\cS) \cap \cT) = \cT$. Hence, 
$f(\cT) = f(\cS \backslash \cK) + \sum\limits_{\ell_{i} \in \cL} f((\cS \cup\{\ell_{1},\dots,\ell_{i}\}) \backslash \cK) - f((\cS \cup\{\ell_{1},\dots,\ell_{i-1}\}) \backslash \cK).$

Given $\ell_{i} \in \cL$, let $\pi_{j} = \ell_{i}$. Then $(\cS \cup \{\ell_{1},\dots,\ell_{i-1}\})\backslash \cK \subseteq \cS \cup \{\pi_{1},\dots,\pi_{j-1}\}$; it follows from the definition of $\dknap[f]$ that
\begin{align*}
&f((\cS \cup\{\ell_{1},\dots,\ell_{i}\}) \backslash \cK) - f((\cS \cup\{\ell_{1},\dots,\ell_{i-1}\}) \backslash \cK) 
\geq \ f(\cS \cup \{\pi_{1},\dots,\pi_{j}\}) - f(\cS \cup \{\pi_{1},\dots,\pi_{j-1}\}) - \dknap[f].
\end{align*}
Therefore,
\begin{align*}
f(\cT) \geq f(\cS \backslash \cK) - |\cL|\dknap[f] + \sum\limits_{\pi_{j} \in \cL} f(\cS \cup \{\pi_{1},\dots,\pi_{j}\}) - f(\cS \cup \{\pi_{1},\dots,\pi_{j-1}\}).
\end{align*}
By the definition of $U_{\pi}(\cS),$ for all $s \in \cS$,
$f(\cS \cup \{\pi_{1},\dots,\pi_{j}\}) - f(\cS \cup \{\pi_{1},\dots,\pi_{j-1}\}) \geq f(\{s\}).$
Because $\cT = (\cS \cup \cL) \backslash \cK$ and $|\cT| \geq |\cS|$, $|\cS \cap \cT| + |\cK|$ $= |\cS|$ $\leq |\cT|$ $= |\cS \cap \cT| + |\cL \cap \cT|$ $= |\cS \cap \cT| + |\cL|$; thus, $|\cK| \leq |\cL|$.
This implies that
\begin{align*}
&f(\cS \backslash \cK) - |\cL|\dknap[f] + \sum\limits_{\pi_{j} \in \cL} f(\cS \cup \{\pi_{1},\dots,\pi_{j}\}) - f(\cS \cup \{\pi_{1},\dots,\pi_{j-1}\})
\geq \ f(\cS \backslash \cK) - |\cL|\dknap[f] + \sum\limits_{s \in \cK} f(\{s\}).
\end{align*}
By the definition of $\dknap[f]$,
$f(\cT) \geq f(\cS \backslash \cK) - (|\cL| + |\cK|)\dknap[f] + \sum\limits_{s \in \cK} f((\cS \backslash \cK) \cup \{s\}) - f(\cS \backslash \cK).$
By the monotonicity of $f$, $f(\cS) \geq f(\cS \backslash \cK \cup\{s\}) \geq f(\cS\backslash \cK)$, for all $s \in \cK$. 
Also, $|\cL| \leq |U_{\pi}(\cS)|$ and $|\cK| \leq |\cS|$. 
Thus, 
$f(\cT) \geq \ f(\cS) - (|\cS| + |U_{\pi}(\cS)|)\dknap[f]
> \ b,$
which follows from the hypothesis. It follows that $f(x(\tilde{\cS})) > b$.

To prove the facet claim, observe that each of the points $x(\cS \backslash \{s\}),$ for all $s \in \cS$ and $x(\cS \cup \{s\} \backslash\{t_{s},u_{s}\})$, for all $s \in U_{\pi}(\cS)$, are $|E_{\pi}(\cS)|$ affinely independent points in $\{x \in X \setbar x_{s} = 0, \forAll s \not\in E_{\pi}(\cS)\}$, and the valid inequality holds with equality for these points. Thus, the valid inequality defines a facet of $\{x \in X \setbar x_{s} = 0, \forAll s \not\in E_{\pi}(\cS)\}$.
$\myQED$

\begin{makeshiftResult}{knapsackFunctionBound}
We have $\dknap[f] \leq p||w||_{1}^{p-1}||w||_{\infty}$. 
\end{makeshiftResult}
\myProof
Let $\cW = \{s \in \Omega \setbar w_{s} > 0\}.$ Let $\ell \in \{0,\dots,|\Omega|-1\}, k \in \{0,\dots,|\Omega|\}$, $\cA, \cB \subseteq \Omega, s \in \Omega, |\cA| = \ell, |\cB| = k$. Suppose that $s \not\in \cW$. Then $(f(\cA \cup \cB \cup \{s\}) - f(\cA \cup \cB)) - (f(\cA \cup \{s\}) - f(\cA)) = u_{s} - u_{s} \leq p||w||_{1}^{1 - p}||w||_{\infty}$. 

Next, suppose $s \in \cW$. Observe that $H$ is Lipschitz continuous on $[0, ||w||_{1}]$, the codomain of $G$: for any $z_{1}, z_{2} \in [0, ||w||_{1}]$, we have $|H(z_{2}) - H(z_{1})| \leq ||H'||_{\infty}|z_{2} - z_{1}| \leq p||w||_{1}^{p-1}|z_{2} - z_{1}|.$ Hence, $f(\cA \cup \cB \cup \{s\}) - f(\cA \cup \cB) \leq u_{s} + p||w||_{1}^{p-1}w_{s} \leq u_{s} + p||w||^{p-1}_{1}||w||_{\infty}.$ Also, $f(\cA \cup \{s\}) - f(\cA) \geq u_{s}$. Hence, $(f(\cA \cup \cB \cup \{s\}) - f(\cA \cup \cB)) - (f(\cA \cup \{s\}) - f(\cA)) \leq p||w||_{1}^{p-1}||w||_{\infty}$, so $\dknap[f] \leq p||w||_{1}^{p-1}||w||_{\infty}$.
$\myQED$

\section{Example of Bounds on the Marginal Violation: The Cooperative Uncapacitated Facility Location Problem} \label{sec:facLocExample}

We present a generalization of the well-known uncapacitated facility location problem (see \cite{Mirchandani} for a detailed overview). We choose uncapacitated facility location as a demonstrative example because of its historical importance (e.g., \cite{Cornuejols77}). In this generalized facility location problem that we consider, the objective function is not submodular in general. We show that the marginal violation metric $D$ (\thref{def:marginalViolation}) can be bounded analytically by exploiting the problem structure and that the objective function's proximity to submodularity is influenced by certain problem parameters.

The objective function of the uncapacitated facility location problem (UFLP) provides an example of a submodular function. An instance of UFLP is defined by $m$ facility locations ($\Omega = \{1,...,m\})$, $n$ clients, demands $b \in \mathbb{R}^{n}_{+}$, fixed costs $w \in \mathbb{R}^{m}_{+}$, and facility-client revenues $v \in \mathbb{R}^{m \times n}$. We consider instances in which $v$ is nonnegative. Additionally, we assume that $w = 0$ so that the firm only assigns facilities to clients based on the variable revenue. We note that \cite{Cornuejols77} consider similar conditions. Let $f: 2^{\Omega} \to \mathbb{R}$ be the objective function of the UFLP with cardinality parameter $K \in \{1,\dots,|\Omega|\}$
\begin{align*}
f(\cS) \coloneqq \begin{cases} \sum\limits_{j = 1}^{n} b_{j}\max\limits_{i \in \cS} v_{ij}, & \mbox{ if } \cS \neq \emptyset\\ 0, & \mbox{ if } \cS = \emptyset. \end{cases}
\quad\quad\quad \text{UFLP:} \ \max\limits_{\cS \subseteq \Omega}\{f(\cS) \ \mathrm{subject \ to} \ |\cS| \leq K\}.
\end{align*}
Here, $\cS$ is a subset of facility locations. Under these conditions, $\myf$ is nonnegative, increasing, and submodular. We consider a generalization of UFLP where the objective function is approximately submodular function. Let $\cS^{2} = \{(p,q) \in \{1,...,m\}^{2},$ for any $\cS \subseteq \Omega$. We introduce a nonnegative reward $u_{pq}$ associated with the simultaneous selection of facilities $p$ and $q$, where $(p,q) \in \Omega^{2}$. We assume that $u_{pp} = 0$ for all $p \in \Omega$. Define $h: 2^{\Omega} \to \mathbb{R}$ as the objective function of the \textit{cooperative uncapacitated facility location problem} (CUFLP) with maximum cardinality parameter $K$.
\begin{align*}
h(\cS) \coloneqq \begin{cases} \sum\limits_{j = 1}^{n} b_{j}\max\limits_{i \in \cS} v_{ij} + \sum\limits_{(p,q) \in \cS^{2}} u_{pq}, & \mbox{ if } \cS \neq \emptyset\\
0 & \mbox{ if } \cS = \emptyset. \end{cases}
\quad\quad\quad \text{CUFLP}: \ \max\limits_{\cS \subseteq \Omega}\{h(\cS) \ \mathrm{subject \ to} \ |\cS| \leq K\}.
\end{align*}

\begin{remark}\thlabel{CUFLP_NPhard}
It is well known that UFLP is NP-hard \citep{Nemhauser83}; thus, CUFLP (which includes UFLP as a special case) is also NP-hard. In addition, the objective function of CUFLP, $h$ is not submodular in general.
\end{remark}

\begin{example}\thlabel{example:CUFLPNotSubmod}
To illustrate the second statement of \thref{CUFLP_NPhard}, consider an instance of the cooperative uncapacitated facility location problem in which $m = 3, n = 1$, and $v_{i1} = 0$, for $i = 1, 2, 3$, $b_{1} = 1,$ $u_{2,3} = 1$, and $u_{pq} = 0$ otherwise. The fixed costs are zero so $\myh$ is increasing. Consider $\cA = \{1\}, \cB = \{2\},$ and $s = \{3\}$. Then,
%
$h(\cA \cup \cB \cup \{s\}) = 1, \quad h(\cA \cup \cB) = 0, 
h(\cA \cup \{s\}) = 0, \quad \text{and} \quad h(\cA) = 0 
\Rightarrow 1 = h(\cA \cup \cB \cup \{s\}) - h(\cA \cup \cB) - h(\cA \cup \{s\}) + h(\cA). $
By \thref{def:Submod}, $\myh$ is not submodular.
\end{example}

Let $\mathrm{supp}(u) \coloneqq \{(p,q) \in \Omega^{2} \setbar u_{pq} > 0\}$.

\begin{proposition}\thlabel{nearSubmodFLP}  
Given an instance of CUFLP, we have $d^{\ell, k}[h]$ $\leq |\mathrm{supp}(u)|\max\{u_{pq} \setbar (p,q)$ $\in \Omega^{2}\}$ for all $\ell \in \{0,...,m - 1\}, k \in \{0,...,m\}.$ Hence, $\dknap[h] \leq |\mathrm{supp}(u)|\max\{u_{pq} \setbar (p,q)$ $\in \Omega^{2}\}$.
\end{proposition}

\begin{proof}
Because $u
$ is nonnegative, $f(\cS) \leq h(\cS)$ for all $\cS \subseteq \Omega$. Further, $f(\cS) \geq h(\cS) - |\mathrm{supp}(u)|\max\{u_{pq} \setbar (p,q) \in \Omega^{2}\}.$ Let $\cA, \cB \subseteq \Omega, s \in \Omega \backslash \cA$, where $|\cA| = \ell \in \{0,...,m - 1\}$ and $|B| = k \in \{0,...,m\}$. 
\begin{align*}
&h(\cA \cup \cB \cup \{s\}) - h(\cA \cup \cB) - h(\cA \cup \{s\}) + h(\cA) \\
\leq \ &h(\cA \cup \cB \cup \{s\}) - f(\cA \cup \cB) - f(\cA \cup \{s\}) + h(\cA)\\
\leq \ &f(\cA \cup \cB \cup \{s\}) - f(\cA \cup \cB) - f(\cA \cup \{s\}) + f(\cA) \\
+ &|\mathrm{supp}(u)|\max\{u_{pq} \setbar (p,q) \in \Omega^{2}\}\\
\leq \ &|\mathrm{supp}(u)|\max\{u_{pq} \setbar (p,q) \in \Omega^{2}\}. 
\end{align*}
It follows that 
$d^{\ell, k}[f] \leq |\mathrm{supp}(u)|\max\{u_{pq} \setbar (p,q) \in \Omega^{2}\}.$ The bound for $\dknap[h]$ follows immediately.
\end{proof}

\thref{nearSubmodFLP} is an example of how one can use the structure of the approximately submodular function in order to derive bounds for the metrics. This allows one to use these bounds immediately as a substitute for the exact value of the metric and avoid exact computation. 

\section{Solution Times for Randomly Generated Approximately Submodular Packing Problem Instances}\label{appx:times}
Table~\ref{tab:InstanceSolveTime} shows solution times for the approximately submodular packing problem instances with and without the addition of the proposed valid inequalities.
\begin{table}[H]
\centering
\scriptsize
\csvautobooktabularcenter[table head = \hline \bfseries No. & \bfseries Time w/o Ineqs (s) & \bfseries Time w/ Ineqs (s) & \bfseries Time Ratio & \bfseries \# of Inequalities Added\\\hline]{TableC1.csv}
\caption{Solution times for randomly generated approximately submodular packing problems.}\label{tab:InstanceSolveTime}
\end{table}

\end{document}